%% file: main.tex
\documentclass{amsart}
\usepackage{graphicx}

\usepackage{fullpage}

\newtheorem{theorem}{Theorem}

\newtheorem{lemma}{Lemma}

\newtheorem{remark}{Remark}

\newtheorem{defn}{Definition}
\newcounter{list_count}
\newenvironment{rlist}{\begin{list}
				{(\roman{list_count})}{\usecounter{list_count}
				\setlength{\rightmargin}{\leftmargin}}}{\end{list}}

\def\integers{{\mathbb Z }}
\def\nnintegers{{\integers_{\ge 0}}}

\def\cO{{\rm O}}
\def\co{{\rm o}}
\def\reals{{\mathbb R }}
\def\naturals{{\mathbb N }}
\def\a{\alpha}

\begin{document}

%              PICTURE STUFF
% Look in the 'Windows' menu for the pictures window
% It's like the Scrapbook -- cut and paste pictures
%

%\def\picture #1 by #2 (#3){
%  \vbox to #2{
%    \hrule width #1 height 0pt depth 0pt
%    \vfill
%    \special{picture #3} % this is the low-level interface
%    }
%  }

%\def\scaledpicture #1 by #2 (#3 scaled #4){{
%  \dimen0=#1 \dimen1=#2
%  \divide\dimen0 by 1000 \multiply\dimen0 by #4
%  \divide\dimen1 by 1000 \multiply\dimen1 by #4
%  \picture \dimen0 by \dimen1 (#3 scaled #4)}
%  }

%inclusion of pictures with Illustrator

\def\picill#1by#2(#3)
{\vbox to #2
{\hrule width #1 height 0pt depth 0pt
\vfill\special{illustration #3}}}
\def\fnum@figure{Figura \thefigure}

%
% Note that you can also say, e.g.,
%  \special{postscript xxx yyy zzz}
% to include PostScript graphics in your documents
%

%\def\annie{\scaledpicture 102pt by 239pt (annie scaled 2000)}
%\def\finder{\picture 260pt by 165pt (screen0 scaled 500)}
%\def\stripes{\picture 2.29in by 1.75in (AWstripes)}
%\def\icon{\picture 7in by 7in (icon)}

%\def\binarypic{\scaledpicture 2.85in by 2.24in (binary scaled 750)}
%\def\binarypic{\scaledpicture 1.85in by 2.24in (binary scaled 750)}
%\def\markovpic{\scaledpicture 3.04in by 3.04in (markov scaled 600)}
%\def\markovpic{\scaledpicture 3.04in by 3.04in (markov scaled 600)}
%\def\ioconepic{\scaledpicture 4.04in by 1.94in (interval1 scaled 750)}
%\def\ioctwopic{\scaledpicture 4.54in by 1.58in (interval2 scaled 750)}
%\def\intervalpic{\scaledpicture 7.17in by 5.29in (intervals scaled
%750)}

%\def\figure2{\scaledpicture 1.67in by 0.54in (fig2 scaled 1000)}

%\parindent=1pc

%Now \TeX\ can talk about the Macintosh
%in pictures, as well as words. (See Figure 1.)
%\smallskip
%\indent\finder

%Figure 1.

%\bye
%
%\\\\\\\\\\\\\\\\\\\\\\\\\\\\\\\\\\\\\
%                                                       text : part 1
%/////////////////////////////////////
%
%\setcounter{page}{1}
%\begin{center}
%{\large \sc
%Dynamical Systems Applied to Asymptotic Geometry.}
%\end{center}

%\begin{center}
%{\it A. A. Pinto and D. Sullivan}\\
%Faculdade de Ci\^encias,
%Universidade do Porto, Portugal \\
%City University of New York, U.S.A.
%\end{center}

\title[Dynamical Systems Applied to Asymptotic Geometry]{Dynamical Systems
Applied to Asymptotic Geometry}
\author{A. A. Pinto}
\address[A. A. Pinto]{Faculdade de Ciencias, Universidade
do Porto\\ 4000 Porto, Portugal.}
%% Note the doubled @@:
\email[A. A. Pinto]{aapinto@fc.up.pt}

%% Second author
\author{D. Sullivan}
\address[D. Sullivan]{Einstein chair, Graduate Center, City University of
New
York,
and SUNY Stonybrook, New York 11794-3651, U.S.A.}
%%\curraddr[]{}% when away from home
%% Note the doubled @@:
%\email[D. A. Rand]{dar@@maths.warwick.ac.uk}
%% \thanks{Research of the second author was supported in part by NSF
%% grant CCR-86-75257 and DARPA Contract N00019-89-J-1988.}

\begin{abstract}
In the paper we discuss two  questions about smooth expanding
dynamical systems on the circle. (i) We characterize the sequences
of asymptotic length ratios which occur for systems with H\"older
continuous derivative. The sequence of asymptotic length ratios are
precisely those given by a positive H\"older continuous function $s$
(solenoid function) on the Cantor set $C$ of $2$-adic integers
satisfying a functional equation called the matching condition. The
functional equation for the $2$-adic integer Cantor set is $$ s
(2x+1)= \frac{s (x)} {s (2x)}
\left( 1+\frac{1}{ s (2x-1)}\right)-1.
$$ We also present a one-to-one correspondence between solenoid
functions and affine classes of $2$-adic quasiperiodic  tilings of
the real line that are fixed points of the 2-amalgamation
operator. (ii) We calculate the precise maximum possible level of
smoothness  for a representative of the system, up to
diffeomorphic conjugacy, in terms of the functions $s$ and
$cr(x)=(1+s(x))/(1+(s(x+1))^{-1})$. For example, in the Lipschitz
structure on $C$ determined by $s$, the maximum smoothness is
$C^{1+\alpha}$ for $0 < \alpha \le 1$ if, and only if, $s$ is
$\alpha$-H\"older continuous. The maximum smoothness is
$C^{2+\alpha}$ for $0 < \alpha \le 1$ if, and only if, $cr$ is
$(1+\alpha)$-H\"older. A curious connection with Mostow type
rigidity is provided by the fact that $s$ must be constant if it
is $\alpha$-H\"older for $\alpha > 1$.
\end{abstract}

\maketitle

%VT
\maketitle
\thispagestyle{empty} 
\input{imsmark}
\SBIMSMark{2004/06}{December 2004}{}

\tableofcontents

\section{Introduction}
\label{jhvbmnk}

One could say that this paper is about the space
$A (2)$ of sequences
$\{ a_1 ,a_2 ,\ldots \}$ of positive real numbers  satisfying
\begin{rlist}
\item
$a_n /a_m$ is exponentially near 1 if $n-m$ is divisible by a high power of
two,
and
\item
$a_3 ,a_5, a_7
,\ldots$ is constructed from $a_1$ and $a_2 ,a_4 ,a_6 ,\ldots$ by the recursion
\begin{equation}
\label{int1}
a_{2n+1} =
\frac{a_n}
{a_{2n}}
\left( 1+
\frac{1}{a_{2n-1}}
\right) - 1.
\end{equation}
\end{rlist}
The only explicit element in $A (2)$ that we know is $\{
1,1,1,\ldots \}$. However, the following theorem shows that $A (2)$
is a dense
  subset of a separable infinite dimensional complex
Banach manifold of \cite{sullivan:MTFAB}.

\begin{theorem}
\label{abab}
The set $A(2)$ is  canonically isomorphic to
\begin{rlist}
\item
[A)] the set of all possible affine structures on the leaves of
the dyadic solenoid $\tilde{S}(2)$ that are transversely H\"older
continuous and invariant by the natural dynamics $\tilde{E}
(2):\tilde{S} (2) \to \tilde{S} (2)$.
\item
[B)] the set of all $C^r$ structures on the circle $S$
invariant by the ``doubling the angle'' expanding dynamics $E(2):S\to S$,
$r>1$.
\item
[C)] the set of all positive H\"older continuous functions $s$ on
the Cantor set $C$ of $2$-adic integers satisfying
$$ s
(2x+1)= \frac{s (x)} {s (2x)}
\left( 1+\frac{1}{ s (2x-1)}\right)-1.
$$
\item
[D)] the set of all affine classes of $2$-adic quasiperiodic
tilings of the real line that are fixed points of the 2-amalgamation
operator.
\item
[E)]  the set of all affine classes of $2$-adic quasiperiodic fixed
grids of the real line.
\end{rlist}
\end{theorem}

See proof of Theorem \ref{abab} in Section \ref{usus}  (in
\cite{sullivan:MTFAB} are studied the uniformly asymptotically
affine (uaa) and the analytic structures on the circle invariant by
the dynamics of $E(2)$
   leaving   the $C^r$ case  for this paper). The
connection between the sequences of $A (2)$ and $C)$ appears from
restricting $s$ in $C)$ to the dense subset of natural numbers in
the Cantor set of $2$-adic integers. The connection between the
sequences of $A)$ with $D)$ and $E)$ follows from the existence of a
dense leaf in the solenoid, with a natural binary grid, which   is
expanded by the dynamics in a manner combinatorially like $x\to 2x$
acting on $\left\{ n/2^k \right\} \subset
\reals$. Then the connection between $A(2)$ and $A)$
follows from using the sequences
$\{ a_0 ,a_1 ,a_2 ,\ldots \}$ to define ratios of
consecutive lengths between integral points of the grid. The
functional equation makes the doubling map look affine between the
integral grid and its double. The $2$-adic  continuity allows the
complete affine structure induced by pullback to impress itself on
the other leaves of the solenoid $\tilde{S} (2)$. The passage from
A) to B) uses the fact that the solenoid $\tilde{S} (2)$ with its
dynamics is the inverse limit system associated to the diagram
\[
\begin{array}{cccccc}
\cdots  \stackrel{E(2)}{\longrightarrow}  &S &
\stackrel{E(2)}{\longrightarrow}& S
& \stackrel{E(2)}{\longrightarrow} & S\\
& ~~~~~~ \uparrow E(2) &
&~~~~~~ \uparrow E(2) & &~~~~~~ \uparrow E(2)\\
\cdots  \stackrel{E(2)}{\longrightarrow}  &S &
\stackrel{E(2)}{\longrightarrow} & S &
\stackrel{E(2)}{\longrightarrow} & S\\
\end{array}
\]
Thus, $\tilde{S} (2)$ projects to $S$, and the affine structures on the
leaves of $\tilde{S} (2)$ determine a {\it canonical family of
solenoidal charts on} $S$ invariant by the dynamics $E(2)$ on $S$.
This canonical family of charts on the circle is compact modulo
affine normalization.  The connection between $B)$ and  $A)$
associates to each $C^{r}$ structure $U$ of the circle $S$ invariant
by $E(2)$  a unique
canonical family of solenoid charts $\mathcal{F}_U$ with the property that
the solenoidal charts are contained  in the structure $U$.
The conection between $B)$ and $C)$ is given by an explicit construction of
a solenoid function $s_U$ using the expanding property of $E(2)$ with
respect to the
$C^{r}$ structure $U$  (see Lemmas \ref{nnnewww} and \ref{nnneeewww}).

In order to state the next theorem, we introduce the following definitions.
 The {\it  ultra-metric $|{\bf u}|_s:C \times C \to \reals^+_0$} is defined as follows.
Let $x=\sum_{m=0}^\infty x_m 2^m  \in C$ and
$y=\sum_{m=0}^\infty y_m 2^m \in C$ be such that
$x_{n} \ldots x_0=y_{n} \ldots y_0$ and $x_{n+1} \ne y_{n+1}$.
For $0 \le i \le n$, let
$A_i = \sum_{m=0}^i x_m 2^m$ and
$E_i = \sum_{m=0}^i  2^m$.
We define
$$|{\bf u}|_s(x,y)=\inf_{0 \le i \le n}
\left\{ 1+  \sum_{j=A_i}^{E_i} \prod_{l=A_i}^j s(l)
+  \sum_{j=0}^{A_i-1}  \prod_{l=j}^{A_i-1}  s(l)
\right\} \ .
$$ 
We present a  geometric interpretation of the ultra-metric  in Section \ref{jrufcp}.
For $\beta > 0$, we say that a function $f:C \to \reals$ is
{\it $\beta$-H\"older}, with respect to the metric $|{\bf u}|=|{\bf u}|_s$, if there
is a constant $d \ge 0$ such that $|f(y)-f(x)| \le d
\left(|{\bf u}|(x,y)\right)^\beta$ for all $x,y \in C$. We say that $f$
is {\it $\beta$-h\"older}, with respect to the metric $|{\bf u}|$, if
there is a continuous function $\epsilon:\reals^+_0 \to
\reals^+_0$, with $\epsilon(0)=0$, such that $|f(y)-f(x)| \le
\epsilon \left(|{\bf u}|(x,y) \right)\left(|{\bf u}|(x,y)\right)^\beta$  for
all $x,y \in C$. By $f$ being {\it Lipschitz} we mean that $f$ is
$1$-H\"older, and by $f$ being {\it lipschitz} we mean that $f$ is
$1$-h\"older.
  Of course on the real line, with respect to the Euclidean metric,
$\beta$-H\"older for $\beta >1$ or lipschitz implies constancy.
We define the  {\it solenoid cross ratio function}
$cr(x):C \to \reals^+$  by 
$cr(x)=(1+s (x))(1+(s(x+1))^{-1})$.

\begin{theorem}
\label{tata}
For every $C^{r}$ structure $U$ of the circle $S$ invariant
by $E(2)$, the overlap maps and  the expanding map
$E(2):S\to S$  attain  its maximum of smoothness with respect to the
canonical family of solenoid charts $\mathcal{F}_U$ contained in $U$.
  Table 1 presents   explicit conditions
in terms of the  corresponding solenoid function
$s=s_U$ which determine the   degree of smoothness of
the overlap homeomorphisms and of $E(2)$
in $\mathcal{F}_U$, and vice-versa.
\bigskip

\begin{tabular}{|c|c|}
\hline
\multicolumn{1}{|p{6.0cm}}
{The regularity of the solenoidal
chart overlap maps and
$E(2):S\to S$.}
&
\multicolumn{1}{|p{5.6cm}|}
{Condition on the functions
$s$  and  $cr$, using the
$|{\bf u}|_s$  ultra-metric on  $C$.}
\\
\hline \hline
have $\a$-H\"older $1^{{\rm st}}$ derivative  &
  $s$ is $\a$-H\"older
\\
$0<\a \le 1$ &
\\
\hline
have $\a$-H\"older $1^{{\rm st}}$ derivative  &
$cr$ is $\a$-H\"older
\\
$0<\a \le 1$ &
\\
\hline
have Lipschitz $1^{{\rm st}}$ derivative &
$s$   is Lipschitz
\\
\hline
have $\a$-H\"older $2^{{\rm nd}}$ derivative  &
$cr$ is $(1+\a)$-H\"older
\\
$0<\a \le 1$ &
\\
\hline have Lipschitz $2^{{\rm nd}}$ derivative & $cr$ is
$2$-H\"older
\\
\hline Affine   & $s$   is lipschitz
\\
\hline
\end{tabular}

\bigskip

\centerline{Table $1$.}

\end{theorem}

See proof of Theorem \ref{tata}  in Section \ref{jrufcp}.
The  scaling and solenoid functions
   give  a   deeper
understanding of the smooth structures  of one dimensional
dynamical systems   (cf. \cite{alby:ss}, \cite{CHJ}, \cite{feib},
\cite{pinto:markov_map}, \cite{sullivan:scalingfn} and \cite{vul})
and also of  two dimensional
dynamical systems (cf. \cite{HR}  and \cite{measures}).

\subsection{Smoothness of diffeomorphisms and ratio distortions of grids}
\label{Int2}

To prove Theorem \ref{tata}, we show some of the relations proposed
in \cite{sullivan:nested} between distinct degrees of smoothness of
a homeomorphism of a real line with distinct bounds of the ratio and
cross ratio distortions of intervals   of a fixed grid that we pass
to describe.

Given $B \ge 1$, $M>1$ and $\Omega:\naturals \to \naturals$,
a $(B,M)$ {\it grid}
$${\mathcal G}_\Omega=\{I_\beta^n \subset I: n \ge 1~{\rm
and}~\beta=1,\ldots,\Omega(n) \}$$
{\it of a closed interval $I$}
is a collection of {\it grid intervals  $I_\beta^n$  at   level   $n$}
with the following properties:
(i)
The grid intervals are closed intervals;
(ii)
For every $n \ge 1$,
the union
$\cup_{\beta=1}^{\Omega(n)} I_\beta^n$
of all   grid  intervals
$I_\beta^n$, at   level   $n$,
is equal to the interval $I$;
(iii)
  For every $n \ge 1$, any two distinct grid intervals  at  level $n$
have disjoint interiors; (iv) For every      $1\le \beta <
\Omega(n)$, the intersection of the grid intervals $I_\beta^n$ and
$I_{\beta+1}^n$ is only an endpoint common to both intervals; (v)
For every $n \ge 1$, the set of all
endpoints of the intervals $I_\beta^n$ at   level  $n$ is
contained in the set of all end points of the intervals
$I_\beta^{n+1}$ at   level  $n+1$; (vi)
  For every    $n \ge 1$ and for every $1\le \beta <
\Omega(n)$, we have
$B^{-1} \le |I_{\beta+1}^n|/|I_\beta^n| \le B$;
(vii)
For every $n \ge 1$ and for
every  $1\le \alpha \le \Omega(n)$,
the   grid interval $I_\alpha^n$
contains at least two grid intervals at level $n+1$,
and contains at most $M$ grid
intervals also at level $n+1$.

Let $h:I \to J$ be a homeomorphism between
two compact  intervals $I$ and $J$ on the real line,
and let ${\mathcal G}_\Omega$ be a grid of $I$. Let $I_\beta$ and
$I_{\beta'}$ be two intervals contained in the real line. The {\it
logarithmic ratio distortion $lrd (I_\beta, I_{\beta'})$}  is given
by $$ lrd (I_\beta, I_{\beta'})   = \log \left( \frac{|I_{\beta}|}
{|I_{\beta'}|} \frac{|h(I_{\beta'})|}{|h(I_{\beta})|} \right) \ . $$
We say that two closed intervals $I_{\beta}$ and $I_{\beta'}$ are
{\it adjacent} if their intersection $I_\beta \cap I_{\beta'}$ is
only an endpoint common to both intervals. Let $I_\beta$,
$I_{\beta'}$ and $I_{\beta''}$ be contained in the real line, such
that $I_\beta$ is adjacent to  $I_{\beta'}$,
  and   $I_{\beta'}$ is adjacent to  $I_{\beta''}$.
The {\it cross ratio} $cr(I_\beta, I_{\beta'}, I_{\beta''})$
is determined by
\begin{eqnarray*}
cr(I_\beta, I_{\beta'}, I_{\beta''})
& = &
\log
\left(
1+
\frac
{|I_{\beta'}|}
{|I_\beta|}
\frac
{|I_\beta|+|I_{\beta'}|+|I_{\beta''}|}
{|I_{\beta''}|}
\right) \ .
\end{eqnarray*}
The  {\it cross ratio distortion} $crd(I_\beta, I_{\beta'},
I_{\beta''})$ is given  by $$ crd (I_\beta, I_{\beta'},
I_{\beta''})= cr(h(I_\beta), h(I_{\beta'}),
h(I_{\beta''}))-cr(I_\beta, I_{\beta'}, I_{\beta''}) \ . $$

\begin{theorem}
\label{th13131313} Let $h:I \to J$ be a homeomorphism between
two compact  intervals $I$ and $J$ on the real line,
and let ${\mathcal G}_\Omega$ be a grid of $I$.
\begin{rlist}
\item
If $h$ has the degree of smoothness presented in
a line of Table $2$, and $dh(x) \ne 0$ for all $x \in I$
(not applicable for quasisymmetric and (uaa) homeomorphisms), then
the logarithmic ratio distortion  satisfy the
bounds presented in the same line with respect to all grid
intervals.
Conversely, if the logarithmic  ratio distortion
satisfies the bounds presented in a line of Table $2$ with respect
to all grid intervals, then $h:I \to J$ has the degree of smoothness
presented in the same line, and $dh(x) \ne 0$ for all $x \in I$ (not
applicable for quasisymmetric and (uaa) homeomorphisms).

\bigskip

\begin{tabular}{|c|c|}
\hline
\multicolumn{1}{|p{5.1cm}|}
{\ \ \ \ \ \ The smoothness of  $h$}
&
\multicolumn{1}{|p{5.1cm}|}
{The order of
$lrd\left(I_\beta^n,I_{\beta+1}^n\right) $}
\\
\hline \hline
 Quasisymmetric & $\cO\left(1\right)$ \\
\hline (uaa) & $\co\left(\left| I_\beta^n \right|\right)\left| I_\beta^n
\right|^{-1}$ 
\\
\hline $C^{1+\alpha}$ & $\cO\left(\left| I_\beta^n \right| ^\alpha \right)$ 
\\
\hline
$C^{1+Lipschitz}$
&
$\cO\left(\left|
I_\beta^n
\right|
\right)$
\\
\hline
Affine
&
$\co\left(\left|
I_\beta^n
\right|
\right)$
\\
\hline
\end{tabular}

\bigskip

\centerline{Table $2$.}

\bigskip

\item
If $h$ has the degree of smoothness presented in
a line of Table $3$, and $dh(x) \ne 0$ for all $x \in I$
(not applicable for quasisymmetric and (uaa) homeomorphisms), then
 the cross ratio distortion satisfy the
bounds presented in the same line with respect to all grid
intervals.
Conversely, if the cross ratio distortion satisfies the bounds presented in a
line of Table $3$ with respect to all grid intervals, then, for
every closed interval $K$     contained  in the interior of $I$, the
homeomorphism $h|K$ restricted to $K$ has the degree of smoothness
presented in the same line, and $dh(x) \ne 0$ for all $x \in I$ (not
applicable for quasisymmetric and (uaa) homeomorphisms).

\bigskip

\begin{tabular}{|c|c|c|}
\hline
\multicolumn{1}{|p{5.1cm}|}
{ \ \ \ \  \   \     The smoothness of  $h$}
&
\multicolumn{1}{|p{5.1cm}|}
{The order of
$crd\left(I_\beta^n,I_{\beta+1}^n,I_{\beta+2}^n\right)$}
\\
\hline \hline 
Quasisymmetric &
$\cO\left(1  \right)$\\
\hline (uaa) &  $\co \left(\left| I_\beta^n \right|\right)\left| I_\beta^n
\right|^{-1}$
\\
\hline $C^{1+\alpha}$  &
$\cO \left(\left| I_\beta^n \right| ^\alpha\right)$
\\
\hline
$C^{2+\alpha}$
&
$\cO\left(\left|
I_\beta^n\right|
^{1+\alpha}\right)$
\\
\hline
$C^{2+Lipschitz}$
&
$\cO\left(\left|
I_\beta^n
\right|
^2\right)$
\\
\hline
\end{tabular}

\bigskip

\centerline{Table $3$.}

\end{rlist}
\end{theorem}

In Section
\ref{section10}, we   present   the  definitions of the degrees of
smoothness presented in Tables $2$ and $3$, and we prove Theorem
\ref{th13131313} in Section \ref{pppttt}. We point out that some of
the difficulties and usefulness of these results come from the fact
that (i) we just compute the bounds of the ratio and cross ratio
distortions
  with respect to a countable set of intervals fixed by a  grid,
and (ii) we do not restrict the grid intervals, at the same level,
to have  necessarily the same lengths. In hyperbolic dynamics, these
grids are naturally determined by Markov partitions.

In \cite{gard2}, \cite{welingtonsebastian} and
\cite{pinto:markov_convergence} other relations
are also  presented     between distinct degrees of smoothness
of a homeomorphim  of the real line with
distinct bounds of ratio and cross ratio distortions of intervals.

\subsection{Interval arithmetics}
Throughout the paper, we use the notation $\phi \leq \cO (\psi
(x))$ to indicate that for all $x$, $|\phi (x)| < c |\psi (x)|$
where $c \geq 1$ is a constant  depending only upon quantities
that are explicitly mentioned. Thus, $\phi (n) < \cO (\mu ^n)$
means that $|\phi (n)| < c\mu ^n$ for some constant $c$ as above.
We also use the notation of interval arithmetic for some
inequalities where:
\begin{rlist} \item if $I$
and $J$ are intervals then $I+J$, $I .  J$ and $I/J$ have the
obvious meaning as intervals, \item if $I=\{ x\}$ then we often
denote $I$ by $x$, and \item $I\pm \epsilon$ denotes the interval
consisting of those $x$ such that $|x-y|<\epsilon$ for some $y\in
I$.
\end{rlist}
Thus $\phi (n) \in 1\pm {\mathcal O}(\nu^n)$ means that there
exists a constant $c>0$ depending only upon explicitly mentioned
quantities such that for all $n\ge 0$, $1-c\nu^n<\phi
(n)<1+c\nu^n$. Similarly, the notation $\phi \leq \co (\psi (x))$
indicates   that for all $x$, $|\phi (x)| < \epsilon( |\psi
(x)|)|\psi (x)|$ where $\epsilon:\reals^+_0 \to \reals^+_0$ is a
continuous function, with $\epsilon(0)=0$,  depending only upon
quantities that are explicitly mentioned.

\section{Expanding dynamics of the circle}
\label{secc2}

In this section we prove a more general version of
Theorems \ref{abab} and \ref{tata} applicable to expanding circle
maps of degree $d$, with $d \geq 2$.

\subsection{$C^{1+ H\ddot older}$ structures $U$ for  the expanding
circle  map $E$} \label{section2} \label{section7}

In this section, we   present the definition of a $C^{1+ H\ddot
older}$  expanding circle  map $E$ with respect to a structure $U$
and give its characterization in terms  of the ratio distortion of
$E$ at small scales with respect to the charts in $U$.

The {\it expanding circle  map }
$E=E(d):S \to S$
with  degree $d \ge 2$ is
given by
$E(z)=z^d$ in complex notation.
Let $p \in S$ be one of the
  fixed points of the
expanding circle map $E$. The {\it  Markov intervals of the
expanding circle map} $E$ are the adjacent closed intervals
 $I_0,\ldots,I_{d-1} $ with non empty interior such that only their
boundaries are contained in the set $\{E^{-1}(p)\}$ of pre-images of
the fixed point $p \in S$. Choose the interval $I_0$ such that
  $I_0 \cap I_{d-1}= \{p\}$.
Let the {\it branch  expanding
circle map} $E_i:I_i \to S$
be   the restriction of
the expanding circle map $E$
to the Markov interval $I_i$,
for all $0 \le i <d$.
Let the interval
$I_{\a_1\ldots\a_n}$ be
$ E_{\a_n}^{-1} \circ
\ldots \circ E_{\a_1}^{-1}
(S).
$
The {\it $n^{\rm th}$-level of the
interval partition
of the expanding circle  map}
$E$
is the set of all
closed intervals
$I_{\a_1\ldots\a_n} \in S$.

A   $C^{1+ H\ddot older}$ diffeomorphism $h:I  \to J$
is a    $C^{1+\epsilon}$ diffeomorphism for some $\epsilon >
0$ (the notion of a quasisymmetric homeomorphism and of a
$C^{1+\epsilon}$ diffeomorphism  $h:I  \to J$ are presented,
respectively,  in sections \ref{section11} and \ref{hhholderrr}.)

\begin{defn}
{\it The  expanding circle map
$E:S \to S$
is   $C^{1+ H\ddot older}$ with respect to a  structure  $U$
on the circle  $S$}
if for every finite cover $U'$ of $U$,
\begin{rlist}
\item
there is an $\epsilon >0$
with the property that
for all  charts
$u:I \to \reals$
and
$v:  J   \to \reals$
contained in $U'$
and for all intervals $K \subset I$
such that $E(K) \subset J$,
the  maps
$v \circ E \circ u^{-1}|u(K)$
are $C^{1+\epsilon}$ and
their $C^{1+\epsilon}$ norms
are bounded away from zero and
infinity;
\item
for every chart $u:I \to    \reals$ contained in $U'$ and for every
map $u_{iso}:I \to \reals$, which is an isometry with respect to the
lengths on the circle $S \subset \reals^2$ determined by the
Euclidean norm on $\reals^2$, the composition $u_{iso} \circ u^{-1}$
is a quasisymmetric homeomorphism.
\end{rlist}
\end{defn}

We note that the above condition (ii) is equivalent to demand that
there are constants $c>0$ and $\nu >1$ such that, for every $n>0$
and every $x \in S$, $|(v \circ E^n \circ u)'(x)| > c \nu^n$, where
$u:I \to    \reals$ and $v:J \to    \reals$ are any two charts in
$U'$ such that $x \in u(I)$ and $E^n \circ u (x) \in K$.

\begin{lemma}
\label{qsq23}
The expanding circle map
$E:S \to S$ is
   $C^{1+ H\ddot older}$ with respect to a  structure  $U$
if, and only if, for  every finite cover  $U'$ of $U$, there are
constants $0 < \mu < 1$  and $b>1$ with the following property: for
all charts $u :J  \to \reals$ and $v :K  \to \reals$ contained in
  $U'$
and for all adjacent intervals $I_{\a_1\ldots\a_n}$ and $I_{\beta_1
\ldots\beta_n }$ at  level $n$ of the interval partition such that
$I_{\a_1\ldots\a_n}, I_{\beta_1 \ldots\beta_n } \subset J$ and
$E(I_{\a_1\ldots\a_n}), E(I_{\beta_1 \ldots\beta_n }) \subset K$, we
have that
\begin{equation}
\label{12eq2} b^{-1} < \frac{|u(I_{\a_1\ldots\a_n})|}{|u(I_{\beta_1
\ldots\beta_n })|} < b ~~~{\rm and}~~~ \left | \log
\frac{|u(I_{\a_1\ldots\a_n})|~ |v(E(I_{\beta_1 \ldots\beta_n }))|}
{|u(I_{\beta_1 \ldots\beta_n })| ~|v(E(I_{\a_1\ldots\a_n}))|} \right
| \le \cO(\mu^n).
\end{equation}
\end{lemma}
\noindent
Lemma \ref{qsq23} follows from Theorem \ref{th13131313}
in   Section \ref{section10}.
%
%

%\begin{lemma}
%\label{qsq25}
By using the Mean Value Theorem we obtain the following
result for a $C^{1+ H\ddot older}$   expanding circle map
$E:S \to S$        with respect to a  structure  $U$.
For every finite cover  $U'$ of $U$,
there is an $\epsilon > 0$, with the property that
for all  charts $u :J  \to \reals$
and $v :K  \to \reals$ contained in
  $U'$
and for all adjacent intervals $I$ and $I'$, such that $I, I'
\subset J$ and $E^n(I), E^n(I') \subset K$, for some $n \ge 1$, we
have
\begin{equation}
\label{12eq3} \left | \log \frac{|u(I)||v(E^n(I'))|}
{|u(I')||v(E^n(I))|}
\right | \le \cO(|v(E^n(I)) \cup v(E^n(I'))|^\epsilon).
\end{equation}
%\end{lemma}

\subsection{Solenoids (\~{E},\~{S})} \label{section5}
\label{section5.1}

In this section, we introduce the notion of a (thca) solenoid
$(\tilde{E},\tilde{S})$  and we prove that a $C^{1+H\ddot older}$
expanding circle map $E$ with respect to a structure $U$
determines a unique (thca) solenoid.

The sequence ${\bf x} = (\ldots,x_3,x_2,x_1,x_0)$
is {\it an  inverse path  of the expanding circle map $E$}
if $E(x_n)= x_{n-1}$, for all $n \ge 1$.
The {\it topological solenoid}
  $\tilde{S}$ consists of
all inverse paths ${\bf x} = (\ldots,x_3,x_2,x_1,x_0)$ of the
expanding circle map $E$ with the product topology. The
topological solenoid is a compact set and is the twist product of
the circle $S$ with the  Cantor set
$\{0,\ldots,d-1\}^{\nnintegers}$. The {\it solenoid map} $\tilde{E}$
is the bijective map defined by
$$
\tilde{E} ({\bf x})=
(\ldots x_0, E(x_0)).
$$
The {\it projection map}
$\pi=\pi_S: \tilde{S} \to S$
is defined by $\pi({\bf x})=x_0$.
A {\it fiber} over $x_0 \in S$ is the set of all
points ${\bf x} \in \tilde{S}$
such that $\pi({\bf x})=x_0$.
A fiber is topologically a
Cantor set $\{0,\ldots,d-1\}^{\nnintegers}$.
A {\it leaf} ${\mathcal  L = L}_{\bf z}$ is the set of all points
${\bf w} \in \tilde{S}$ path connected to the point ${\bf z} \in
\tilde{S}$. A {\it local leaf} ${\mathcal L'}$ is a path connected
subset of a leaf. The {\it monodromy map} $\tilde{M}:\tilde{S} \to
\tilde{S}$ is defined such that the local leaf starting on $\bf x$
and ending on $\tilde{M}({\bf x})$ after being projected by $\pi$ is
an
anti-clockwise arc
starting
on $x_0$, going around the circle once,
and ending   on the point $x_0$.
Since the orbit of any  point ${\bf x} \in \tilde{S}$
under $\tilde{M}$ is dense on its fiber
(see   Lemma
\ref{fdswqewdsfs} in Section \ref{mdadic}),
we get that all leaves $ \mathcal  L$
of the solenoid
$\tilde{S}$ are dense.

\begin{defn}
{\it The solenoid
$(\tilde{E},\tilde{S})$
  is  transversely H\"older continuous affine
(thca)} if  (i) every leaf $\mathcal L$ has an affine structure;
(ii) the solenoid map $\tilde{E}$ preserves the affine structure
on the leaves;
and (iii) the ratio between adjacent leaves
determined by their affine structure
changes H\"older continuously along transversals.
\end{defn}

We say that  {\it $({\bf x,y,z})$  is a   triple}, if the points
${\bf x}$,  ${\bf y}$ and ${\bf z}$ are distinct and are contained
in the same leaf $\mathcal L$ of $\tilde{S}$. Let $T$ be the set of
all triples
$({\bf x,y,z})$.
A function $r:T\to \reals^+$ is {\it invariant by the action of the
solenoid map $\tilde{E}$} if, and only if, for all triples $({\bf
x},{\bf y},{\bf z}) \in T$, we have
$r({\bf x},{\bf y},{\bf z})=
r(\tilde{E}({\bf x}),\tilde{E}({\bf y}),\tilde{E}({\bf z}))$.
A function $r:T\to \reals^+$
{\it varies H\"older continuously
along fibers}
if and only if there are constants $c>0$ and $0 < \mu <1$
with the property  that  for all triples
$({\bf x},{\bf y},{\bf z}),
({\bf x'},{\bf y'},{\bf z'}) \in T$, such that $x_n=x_n'$,
$y_n=y_n'$ and $z_n=z_n'$, we have
$$\left| \log (r({\bf x},{\bf y},{\bf z}) )
- \log (r({\bf x'},{\bf y'},{\bf z'}) )
\right|\le \cO(\mu^n).
$$

\begin{figure}
\centerline{
\includegraphics[width=178pt,height=72pt]{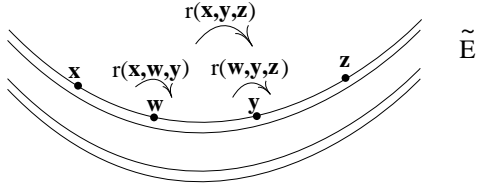}}
% \picill 178pt by 72pt (car06) }
\caption{The matching condition for the leaf ratio function $r$.}
\label{car06}
\end{figure}

\begin{defn}
A {\it   H\"older leaf ratio function} $r:T \to \reals^+$ is a
continuous function invariant by the action of the solenoid map
$\tilde{E}$     that is H\"older continuously along fibers, and
satisfies the following {\it matching condition} (see Figure
\ref{car06}): for all triples $({\bf x}, {\bf w}, {\bf y}), ({\bf
w}, {\bf y}, {\bf z}) \in T$, $$ r({\bf x}, {\bf y}, {\bf z})=
\frac {r({\bf x}, {\bf w}, {\bf y}) r({\bf w}, {\bf y}, {\bf z})}
{1+r({\bf x}, {\bf w}, {\bf y})}. $$
\end{defn}

\begin{lemma}
\label{lasda}
There is a one-to-one
correspondence
between (thca)
solenoids
$(\tilde{E},\tilde{S})$ and
H\"older leaf
ratio functions
$r:T \to \reals^+$.
\end{lemma}

\noindent {\bf Proof:} The affine structures on the leaves of the
(thca) solenoid ${\tilde S}$ determine a function $r:T\to
\reals^+$  that varies continuously along leaves, and satisfies
the matching condition. The converse is also true. Moreover, (i)
the solenoid map $\tilde{S}$ preserves the affine structure on the
leaves if and only if the function $r:T\to \reals^+$ is  invariant
by the action of the solenoid map $\tilde{E}$ and (ii) the ratio
between adjacent leaves determined by their affine structure
changes H\"older continuously along transversals if and only if
the function $r:T\to \reals^+$ varies H\"older continuously along
fibers. \qed

\begin{lemma}
\label{nnnewww} A $C^{1+ H\ddot older}$   expanding circle map
  $E:S \to S$    with respect to a  structure  $U$
  generates a
  H\"older leaf ratio function
$r_U:T \to \reals^+$.
\end{lemma}

\noindent {\bf Proof:} Let $U'$ be a finite cover of $U$. For every
triple $({\bf x},{\bf y},{\bf z}) \in T$ and every $n$ large enough,
let $u_n:J_n \to \reals$ be a chart contained in $U'$ such that
$x_n, y_n ,z_n \in J_n$. Using (\ref{12eq3}),  $r_U({\bf x},{\bf
y},{\bf z})$ is well-defined by
$$
r_U ({\bf x},{\bf y},{\bf z})=
\lim_{n \to \infty}
\frac{|u_n(y_n)-u_n(z_n)|}
{|u_n(x_n)-u_n(y_n)|}.
$$
By construction, $r_U$ is invariant by the dynamics of the solenoid
map and satisfies the matching condition. Again, using
(\ref{12eq3}), we obtain that $r_U$ is a continuous function varying
H\"older continuously along transversals. Hence, $r_U$ is a leaf
ratio function. \qed

\subsection{Solenoid functions $s:C \to \reals^+$}
\label{mdadic}

In this section, we will introduce the notion of a solenoid function
whose domain is a fiber of the solenoid.
We will show that a H\"older leaf ratio function determines
a H\"older solenoid function and that a  H\"older solenoid function
determines an element in the set of sequences $A(d)$.

Let $\sum_{i=-\infty}^\infty a_i d^i$
be a $d$-{\it adic number}.
The $d$-adic numbers
$$
\sum_{i=-\infty}^{n-1} (d-1) d^i
+\sum_{i=n}^\infty a_i d^i
~~~~~{\rm and}~~~~~
(a_n+1) d^n +
\sum_{i=n+1}^\infty a_i d^i
$$
such that  $a_n+1 < d$
are {\it $d$-adic equivalent}.
{\it The $d$-adic set $\tilde{\Omega}$}  is
  the topological Cantor set
$\{0,\ldots, d-1\}^\integers$
of  all $d$-adic numbers
modulo the above
$d$-adic equivalence.
The {\it product  map}
$d\times:\tilde{\Omega} \to \tilde{\Omega}$
is the
multiplication by
$d$ of the $d$-adic numbers.
The {\it add $1$ map}
$1+:\tilde{\Omega} \to \tilde{\Omega}$
is the sum of $1$ to  the $d$-adic numbers.
%and the {\it add $d$ map}
%$d+:\tilde{\Omega} \to \tilde{\Omega}$
%is   $d$ compositions
%of the add $1$ map.

Let the map
$\tilde{\omega}:\tilde{\Omega}  \to \tilde{S}$
be the homeomorphism between the
$d$-adic set $\tilde{\Omega}$ and
the solenoid $\tilde{S}$
defined as follows:
$\tilde{\omega}(\sum_{i=-\infty}^\infty a_i d^i) =
{\bf x}=(\ldots,x_1,x_0) \in \tilde{S}$,
where $x_n= \cap_{i=1}^\infty
E^{-1}_{a_{n-1}} \circ \ldots \circ  E^{-1}_{a_{n-i}} (I _{a_{n-(i+1)}})$
for all $n \ge 0$ (recall that $I_{a_{n-(i+1)}}$ is  a Markov
interval of the expanding circle map $E$). Hence $x_n \in I_{a_n}$
for all $n \ge 0$. By construction, the map
$\tilde{\omega}:\tilde{\Omega} \to \tilde{S}$ conjugates the product
map $d\times:\tilde{\Omega} \to \tilde{\Omega}$ with  the solenoid
map $\tilde{E}:\tilde{S} \to \tilde{S}$, and  conjugates
the   add $1$
map $1+:\tilde{\Omega} \to \tilde{\Omega}$ with the monodromy map
$\tilde{M}:\tilde{S} \to \tilde{S}$.

\begin{lemma}
\label{fdswqewdsfs}
Every orbit of the monodromy map
is dense on its fiber.
\end{lemma}

\noindent {\bf Proof:} Since the add $1$ map $1+:\tilde{\Omega} \to
\tilde{\Omega}$ is dense on the image $\tilde{\omega}^{-1}(F)$ of
every fiber $F$ of the solenoid  $\tilde{S}$, the lemma follows.
\qed

\bigskip

Let  $\Omega$ be
  the topological Cantor set
$\{0,\ldots, d-1\}^{\integers_{\le  0}}$
corresponding to
all $d$-adic numbers
of the form $\sum_{i=-\infty}^{-1} a_i d^i$ modulo the $d$-adic
equivalence. {\it The projection map $\pi_\Omega:\tilde{\Omega} \to
\Omega$} is defined by $\pi_\Omega \left(\sum_{i=-\infty}^\infty a_i
d^i\right)= \sum_{i=-\infty}^{-1} a_i d^i$. The map $\omega : \Omega
\to S$ is defined by $\omega(\sum_{i=-\infty}^{-1} a_i d^i)=
\cap_{i=1}^\infty
E^{-1}_{a_{-1}} \circ \ldots \circ  E^{-1}_{a_{-i}} (I _{a_{-(i+1)}})$.
By construction,
$$\omega\circ  \pi_\Omega
\left(\sum_{i=-\infty}^\infty a_i d^i\right)=
\pi_S \circ   \tilde{\omega}
\left(\sum_{i=-\infty}^\infty a_i d^i\right),$$
for all
$\sum_{i=-\infty}^\infty a_i d^i \in
\tilde{\Omega}$.

{\it The set $C$}  is
  the topological Cantor set
$\{0,\ldots, d-1\}^{\integers_{\ge 0}}$
corresponding to
all $d$-adic integers
of the form $\sum_{i=0}^\infty a_i d^i$.

\begin{figure}
\centerline{
\includegraphics[width=233pt,height=60pt]{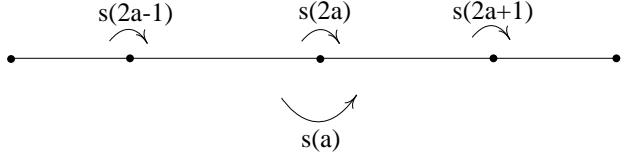}}
%width=0.32\textwidth para 32% da largura (car04)
\caption{The matching condition for the solenoid function ($d=2$).}
\label{car04}
\end{figure}

\begin{defn}
The {\it solenoid function}
$s:C \to \reals^+$
is a continuous function satisfying the following {\it matching
condition} (see Figure \ref{car04}), for all $a \in C$:
\begin{equation}
\label{pdh} s(a) = \frac { \prod_{i=1}^{d-1} s(da-i) \left (
\sum_{j=0}^{d-1} \prod_{l=0}^j s(da+l) \right)} {1+ \sum_{j=1}^{d-1}
\prod_{l=j}^{d-1} s(da-l)} \ .
\end{equation}
\end{defn}

\begin{lemma}
\label{nnneeewww}
The H\"older leaf
ratio function
$r:T \to \reals^+$
determines a H\"older
solenoid function
$s_r:C \to \reals^+$.
\end{lemma}

\noindent {\bf Proof:} For all $\sum_{i=0}^\infty a_i d^i \in C$,
we define
$$
s_r\left(\sum_{i=0}^\infty a_i d^i\right)
=
r \left( \tilde{\omega}
\left(\sum_{i=0}^\infty a_i d^i -1 \right),
\tilde{\omega}
\left(\sum_{i=0}^\infty a_i d^i \right),
\tilde{\omega}
\left(\sum_{i=0}^\infty a_i d^i + 1\right)
\right).$$
The matching condition
and the H\"older continuity of the leaf
ratio function
$r:T \to \reals^+$
imply the
matching condition
and the
H\"older continuity of the
solenoid function
$s_r:C \to \reals^+$,
respectively.
\qed

\begin{lemma}
\label{addddp}
There is a one-to-one
correspondence between H\"older
solenoid functions
$s:C \to \reals^+$
and
sequences
$\{r_1 ,r_2 , r_3, \ldots \} \in  A (d)$
of positive real numbers
with the following properties:
\begin{rlist}
\item
$r_n /r_m \le \cO(\mu^{i})$
if $n-m$ is divisible
by  $d^i$, where  $0 < \mu < 1$;
\item
$r_1 ,r_2, \ldots$
satisfies
\begin{equation}
\label{match} r_a = \frac { \prod_{i=1}^{d-1} r_{da-i} \left (
\sum_{j=0}^{d-1} \prod_{l=0}^j r_{da+l} \right)} {1+
\sum_{j=1}^{d-1} \prod_{l=j}^{d-1} r_{da-l}} \ .
\end{equation}
\end{rlist}
\end{lemma}

A   geometric interpretation
of the sequences contained in
the set $A(d)$ is given by the $d$-adic quasiperiodic tilings and
grids of the real line defined in Section \ref{section10.1}, below.

\bigskip
\noindent {\bf Proof:} Given a H\"older solenoid function $s:C \to
\reals^+$, for all $i=\sum_{j=0}^{k} a_j d^j \ge 0$, we
define $r_i$
by
$$
r_i=s\left(\sum_{j=0}^{k} a_j d^j \right).
$$
The matching condition
of the solenoid function
$s:C \to \reals^+$
implies that
the ratios $r_1,r_2, \ldots$
satisfy  (\ref{match}).
The H\"older continuity
of the solenoid function $s:C \to \reals^+$
implies condition (i).
Conversely,
for every $d$-adic integer
$a=\sum_{i=0}^\infty  a_i d^i \in C$,
let  $\underline{a}_n \in \nnintegers$ be equal to
$\sum_{i=0}^n  a_i d^i$.
Define the value $s(a)$ by
$$
s(a)=\lim_{n \to \infty}
r_{\underline{a}_n}.
$$
Using condition (i) the above limit is well
defined and the function
$s:C \to \reals^+$ is H\"older continuous.
Using condition (ii) and the continuity of $s$
we obtain that  the function
$s$ satisfies the
matching condition. \qed

\subsection{d-Adic quasiperiodic      tilings and grids  and
amalgamation operators} \label{section10.1}

In this section, we introduce  $d$-adic quasiperiodic  tilings of the
real line that are fixed points of the d-amalgamation operator and
$d$-adic quasiperiodic fixed grids of the real line, and we show that
their affine classes are in one-to-one correspondence with (thca)
solenoids.

A {\it tiling} ${\mathcal T}=\{I_\beta \subset \reals: \beta \in
\integers \}$ {\it of the real line} is a collection of {\it tiling
intervals $I_\beta$} with the following properties: (i) The tiling
intervals are closed intervals; (ii) The union $\cup_{\beta\in
\integers} I_\beta$ of all   tiling intervals $I_\beta$ is  equal to
the real line; (iii) any two distinct tiling intervals  have
disjoint interiors; (iv)   For every      $\beta \in \integers$, the
intersection of the tiling intervals $I_\beta$ and $I_{\beta+1}$ is
only an endpoint common to both intervals; (v)
There is $B \ge 1$, such that for every      $\beta \in \integers$, we have
$ B^{-1} \le  |I_{\beta+1}| / |I_\beta|  \le B$. We say that the
tilings  ${\mathcal T}_1=\{I_\beta \subset \reals: \beta \in
\integers \}$ and ${\mathcal T}_2=\{J_\beta \subset \reals: \beta
\in \integers \}$ of the real line are in the same affine class if
there is an affine map $h:\reals \to \reals$ such that
$h(I_\beta)=J_\beta$ for every $\beta \in \integers$.
The {\it tiling sequence} $\underline r=(r_m)_{m \in
\integers}$ is given by $r_m = |I_{m +1}|/|I_m|$. We note that a
tiling sequence $\underline r$ determines an affine class of tilings
$\mathcal T$ and vice-versa. We say that a tiling sequence is {\it
$d$-adic quasiperiodic} if there is $0< \mu <1$ such that  $|r_j -
r_k| \le \cO(\mu^i)$, when $(j-k)/d^i$ is an integer. Let  $\bf T$
denote the set of all tiling
sequences. The {\it d-amalgamation operator} $A_d:\bf T  \to \bf T$ is
defined by  $A_d(\underline r)= \underline s$, where $$ s_i=
r_{d(i-1)+1,di} \frac {1+\sum_{m=di+1}^{d(i+1)-1} r_{di+1,m}}
{1+\sum_{m=d(i-1)+1}^{di-1} r_{d(i-1)+1,m}}, $$ for all $i \in
\integers$.

\begin{defn}
A $d$-adic quasiperiodic    tiling ${\mathcal T}$ of the real line is
a  tiling such that the corresponding tiling sequence $\underline r$
is $d$-adic quasiperiodic. A  tiling ${\mathcal T}$ of the real line
is a fixed point  of the $d$-amalgamation operator
if  the corresponding tiling sequence is a fixed point of the
  d-amalgamation operator, i.e. $A_d(\underline r)=\underline r$.
\end{defn}

\begin{remark}
The  tiling sequence  $\underline r=(r_m)_{m \in \integers}$ of a
$d$-adic quasiperiodic    tiling of the real line that is a fixed
point of the d-amalgamation operator determines a sequence
$r_1,r_2,\ldots$    in $A(d)$.
\end{remark}

A {\it d-grid} ${\mathcal G}$ of the real line is a collection
  of intervals $I_\beta^n$    satisfying properties (i)
to (vii) of a $(B,d)$-grid   ${\mathcal G}_\Omega$, for some $B \ge
1$, such that
  every interval $I_\beta^n$ is the union of
$d$ grid intervals at level $n+1$, and $\Omega(n)=\infty$. We note
that every level $n$ of a grid forms a tiling of the real line. We
say that the grids ${\mathcal G}_1=\{I_\beta^n\}$ and ${\mathcal
G}_2=\{J_\beta^n\}$ of the real line are in the same affine class if
there is an affine map $h:\reals \to \reals$ such that
$h(I_\beta^n)=J_\beta^n$ for every $\beta \in \integers$ and every
$n \in \naturals$. The {\it d-grid sequence} $\ldots \underline r^2
\underline r^1$ is given by $\underline r^n =(r_m^n)_{m \in
\integers}$ where $r_m^n = |I_{m +1}^n|/|I_m^n|$. The following
remark gives a geometric interpretation of the $d$-amalgamation
operator.

\begin{remark}  \begin{rlist} \item If $\ldots
\underline r^2 \underline r^1$ is a d-grid sequence then
$A_d(\underline r^{n+1})= \underline r^n$ for every $n \ge 1$.
\item If $\ldots
\underline r^2 \underline r^1$  is a sequence such that
$A_d(\underline r^{n+1})= \underline r^n$ then the sequence
determines an affine class of d-grids.
\end{rlist}
\end{remark}

\begin{defn}
A $d$-adic quasiperiodic fixed  grid ${\mathcal G}$ of the real line
is a  $d$-grid of the real line such that the corresponding grid
sequence $\ldots \underline r^2 \underline r^1$ is constant, i.e.
$\underline r^1=\underline r^n$ for every $n \ge 1$, and $\underline
r^1$
is $d$-adic quasiperiodic.
\end{defn}

Hence, all the levels of a $d$-adic quasiperiodic fixed grid
${\mathcal G}$ of the real line determine the same $d$-adic
quasiperiodic  tiling of the real line, up to affine equivalence,
that is a fixed point of the $d$-amalgamation operator.

\begin{lemma}
\label{lassssspq} There is a one-to-one correspondence between (i)
(thca) solenoids ; (ii) affine classes of
$d$-adic quasiperiodic    tilings  of the real line that
are fixed points of the d-amalgamation operator; (iii) affine
classes of $d$-adic quasiperiodic fixed   grids   of the
real line.
\end{lemma}

\begin{figure}
\centerline{
\includegraphics[width=348pt,height=145pt]{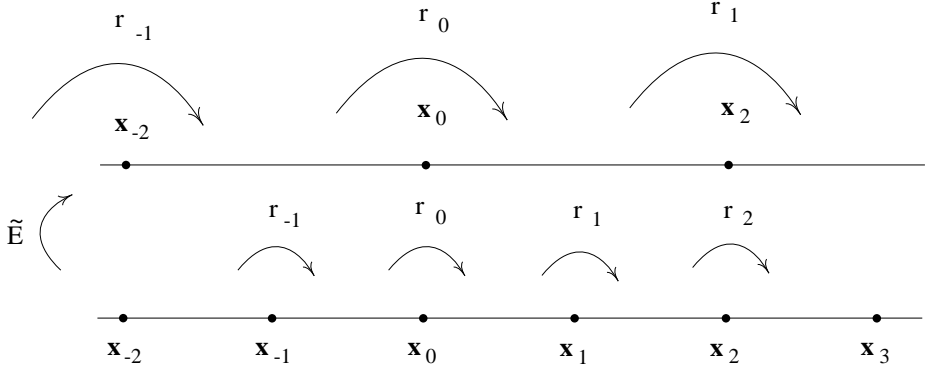}}
% \picill 359pt by 150pt (car12) }
\caption{The leaf
$\mathcal L$ fixed by the solenoid map $\tilde{E}$.} \label{car12}
\end{figure}

\noindent {\bf Proof:} By construction, there is
a one-to-one correspondence between
(ii) affine classes of $d$-adic quasi\-pe\-riodic
tilings    of the real line that are fixed points of the
d-amalgamation operator  and  (iii) affine classes of $d$-adic
quasi\-pe\-riodic fixed grids   of the real line.
Let us prove that a (thca) solenoid determines canonically an affine
class of  $d$-adic quasiperiodic   tilings of the real line that are
  fixed points  of the d-amalgamation operator.  Let $\mathcal L$ be
a leaf of the (thca) solenoid $(\tilde{E},\tilde{S})$  containing
a fixed point ${\bf x}_0$ of the solenoid map $\tilde{E}$. The
leaf $\mathcal L$ is marked by the points $\ldots, {\bf x}_{-1},
{\bf x}_0, {\bf x}_1, \ldots$ that project on the same point of
the circle as the fixed point ${\bf x}_0$, and such that there is
a local leaf ${\mathcal L}_m$ with extreme points ${\bf x}_m$ and
${\bf x}_{m+1}$ with the property that  ${\mathcal L}_m$ does  not
contain any other point ${\bf x}_j$ for $m \ne j \ne m+1$. The
affine structure on the leaf ${\mathcal L}$ determines the ratios
$r_m=r({\bf x}_{m-1},{\bf x}_m,{\bf x}_{m+1})$ of the leaf ratio
function $r:T \to \reals^+$, for all $m \in \integers$. Since the
solenoid map $\tilde{E}$ is affine and $\tilde{E}({\mathcal
L})={\mathcal L}$, the sequence of ratios $\underline r = (r_m)_{m
\in \integers}$ is fixed by the amalgamation operator $A_d$ (see
Figure \ref{car12}). The H\"older transversality    of the
solenoid $(\tilde{E},\tilde{S})$ implies that the sequence
$\underline r$ is $d$-adic quasiperiodic. Therefore, the sequence
$\underline r$ determines an affine class of $d$-adic
quasiperiodic     tilings of the real line that are fixed point of
the d-amalgamation operator, and so an affine class of $d$-adic
quasiperiodic fixed grids of the real line.
Conversely,  an affine class of $d$-adic quasiperiodic fixed grids
of the real line determines uniquely the affine structure of a
leaf $\mathcal L$ that is fixed by the solenoid map $\tilde{E}$.
Since the {\it grid sequence} $\ldots \underline r^2 \underline
r^1$ is a fixed point of the amalgamation operator, i.e.
$A_d(\underline r^n)=\underline r^{n-1}$, the solenoid map
$\tilde{E}$ is affine on the leaf $\mathcal L$. By density of the
leaf $\mathcal L$ on the solenoid $\tilde{S}$ and since the grid
$g_d$ is $d$-adic quasiperiodic, the affine structure of the leaf
$\mathcal L$ extends to an affine structure
transversely H\"older  continuous
on the solenoid  $\tilde{S}$ such that the solenoid map $\tilde{E}$
leaves the affine structure  invariant. \qed

\subsection{Solenoidal charts for the $C^{1+H\ddot older}$
expanding circle map $E$} \label{CCC}

In this section, we introduce the solenoidal charts which  will
determine a canonical structure for the expanding circle map.

\begin{figure}
\centerline{
\includegraphics[width=155pt,height=126pt]{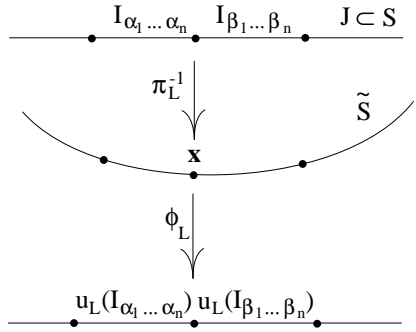}}
% \picill 161pt by 127pt (car08) }
\caption{The solenoidal chart.} \label{car08}
\end{figure}

\begin{defn}
Let ${\mathcal L}$ be a local leaf with an affine
structure and $\pi_{\mathcal L}=\pi_S|{\mathcal L}$ the
homeomorphic projection of  ${\mathcal L}$ onto an interval
$J_{\mathcal L}$ of the circle $S$. Let $\phi_{\mathcal
L}:{\mathcal L} \to \reals$ be a map  preserving the affine
structure of the leaf ${\mathcal L}$.
A {\it solenoidal chart $u_{\mathcal L}:J_{\mathcal L} \to \reals$ on
the circle $S$} is defined by $u_{\mathcal L}= \phi_{\mathcal L}
\circ \pi_{\mathcal L}^{-1}$ (see Figure \ref{car08}).
\end{defn}

\begin{lemma}
\label{rrr} The  solenoidal charts determined by a (thca) solenoid
$(\tilde{E},\tilde{S})$ produce  a canonical structure $U$ such
that
the expanding
circle map $E$ is $C^{1+H\ddot older}$.
\end{lemma}

\noindent {\bf Proof:} Let $U'$ be a finite cover consisting of
solenoidal charts. Let $I_{\a_1\ldots\a_n}$ and
$I_{\beta_1\ldots\beta_n}$ be   adjacent intervals at level $n$ of
the interval partition and $u_{\mathcal L}:J \to \reals$ and
$v_{\mathcal L'}:K \to \reals$ solenoidal charts such that
$I_{\a_1\ldots\a_n}, I_{\beta_1\ldots\beta_n} \subset J$ and
$I_{\a_2\ldots\a_n}, I_{\beta_2\ldots\beta_n} \subset K$. Let ${\bf
x}$, ${\bf y}$ and ${\bf z}$
be the points contained in ${\mathcal L}$
such that $\pi({\bf x})$ and $\pi({\bf y})$ are the endpoints of
$I_{\a_1\ldots\a_n}$, and  $\pi({\bf y})$ and $\pi({\bf z})$ are the
endpoints of $I_{\beta_1\ldots\beta_n}$. Let ${\bf x}'$, ${\bf y}'$
and ${\bf z}'$
be the points contained in ${\mathcal L'}$
such that $\pi({\bf x}')$ and $\pi({\bf y}')$ are the endpoints of
$I_{\a_2\ldots\a_n}$, and  $\pi({\bf y}')$ and $\pi({\bf z}')$ are
the endpoints of $I_{\beta_2\ldots \beta_n}$ (see Figure \ref{car08}). By
Lemma  \ref{lasda},
  the (thca) solenoid determines a leaf ratio
function $r:T \to \reals^+$ such that
\begin{equation}
\label{afk1}
\frac  {|u_{\mathcal
L}(I_{\beta_1\ldots\beta_n})|} {|u_{\mathcal
L}(I_{\a_1\ldots\a_n})|}\frac {|v_{\mathcal
L'}(I_{\a_2\ldots\a_n})|} {|v_{\mathcal
L'}(I_{\beta_2\ldots\beta_n})|} =\frac {r({\bf x},{\bf y},{\bf z})}
{r({\bf x}',{\bf y}',{\bf z}')} \ .
\end{equation}
By Lemma   \ref{nnneeewww}, using that $\tilde{E}$ is affine on
leaves, the leaf ratio function $r:T \to \reals^+$ determines  a
solenoid function $s_r:C \to \reals^+$ such that
\begin{equation}
\label{afk111}
  \frac {r({\bf x},{\bf
y},{\bf z})} {r({\bf x}',{\bf y}',{\bf z}')} = \frac
{s\left(\tilde{\omega}^{-1}(\tilde{E}^n({\bf y}))\right)}
{s\left(\tilde{\omega}^{-1}(\tilde{E}^{n-1}({\bf y'}))\right)} \ .
\end{equation}
By H\"older continuity of the solenoid function,
\begin{equation}
\label{afk2} \left| \log \frac
{s\left(\tilde{\omega}^{-1}(\tilde{E}^n({\bf x}))\right)}
{s\left(\tilde{\omega}^{-1}(\tilde{E}^{n-1}({\bf y}))\right)}
\right| \le \cO(\mu^n),
\end{equation}
for some $0 < \mu < 1$. Putting  (\ref{afk1}), (\ref{afk111}) and
(\ref{afk2}) together, and using that $C$ is compact, we obtain that
\begin{equation}
\label{fffeee} b^{-1} < \frac{|u_{\mathcal L}(I_{\a_1\ldots\a_n})|}
{|u_{\mathcal L}(I_{\beta_1 \ldots\beta_n })|} < b ~~~{\rm and}~~~
\left | \log \frac{|u_{\mathcal L} (I_{\a_1\ldots\a_n})|}
{|u_{\mathcal L}(I_{\beta_1 \ldots\beta_n })|} \frac{ |v_{\mathcal
L'}(I_{\beta_2 \ldots\beta_n })|} {|v_{\mathcal
L'}(I_{\a_2\ldots\a_n})|} \right | \le \cO(\mu^n)
\end{equation}
for some $b \ge 1$. Hence, by Lemma \ref{qsq23}, the expanding
circle map $E$ is $C^{1+H\ddot older}$ with respect to the structure
$U$ produced by the solenoidal charts. \qed

\begin{lemma}
\label{rrr555} The H\"older solenoid function $s:C \to \reals^+$
determines a set of
  solenoidal charts
which produce  a structure $U$ such that
the expanding
circle map $E$ is $C^{1+H\ddot older}$.
\end{lemma}

\noindent {\bf Proof:} For every triple $({\bf x,y,z})$ such that
there are $n \in \integers$ and  $a \in C$ with the property that
$$(\tilde{E}^n({\bf x}),\tilde{E}^n({\bf y}),\tilde{E}^n({\bf z}))
= (\tilde{\omega}(a-1),\tilde{\omega}(a),\tilde{\omega}(a+1))$$ we
define  $r({\bf x,y,z})$ equal to $s(a)$. Hence, the ratios $r$ are
invariant under the solenoid map $\tilde{E}$. Since the solenoid
function satisfies the matching condition, the above  ratios $r$
determine an affine structure on the leaves of the solenoid. By
construction, the solenoidal charts determined by this affine
structure on the leaves satisfy   (\ref{fffeee}), and so by Lemma
\ref{qsq23}, the expanding circle map $E$ is $C^{1+H\ddot older}$
with respect to the structure $U$ produced by the solenoidal charts.
\qed

\subsection{Smooth properties of  solenoidal charts}
\label{section9}

We will prove that the solenoidal charts maximize the smoothness of
the expanding circle map with respect to all charts in the same
$C^{1+H\ddot older}$ structure.

Let $U$ be a $C^{1+H\ddot older}$ structure for the expanding circle
map $E$. By Lemmas \ref{lasda} and \ref{nnnewww}, the structure $U$
determines a (thca) solenoid $(\tilde{E},\tilde{S})_U$.

\begin{lemma}
\label{be}
Let $U$ be a $C^{1+H\ddot older}$ structure for the expanding circle
map $E$, and let $V$ be the
set of all  solenoidal charts  determined by the
(thca) solenoid $(\tilde{E},\tilde{S})_U$.
Then, the set $V$ is contained in $U$  and
the degree of smoothness of
the expanding circle map $E$
when measured in terms of a cover $U'$ of $U$ attains its maximum
when $U' \subset V$.
\end{lemma}

\begin{figure}
\centerline{
\includegraphics[width=177pt,height=107pt]{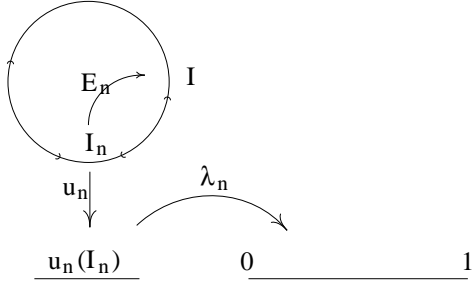}}
% \picill 190pt by 109pt (car09) }
\caption{The construction of the solenoidal charts from
the $C^{1+H\ddot older}$ structure $U$.} \label{car09}
\end{figure}

\bigskip
\noindent {\bf Proof:} Let the  expanding circle map $E:S \to S$ be
$C^r$, for some $r>1$, with  respect to a finite  cover $U'$ of the
structure $U$. We shall prove that the solenoidal charts
$v_{\mathcal L}:I \to \reals$
are $C^r$ compatible with the charts
contained in $U'$,
  proving the theorem.
Let ${\mathcal L}$ be a local leaf that projects by $\pi_{\mathcal
L}=\pi_S|{\mathcal L}$ homeomorphically on an interval $I$
contained in the domain $J$ of a chart $u:J \to \reals$ of $U'$.
For $n$ large enough, let $u_n:J_n \to  \reals$ be a chart in $U'$
such that $I_n=\pi_S(\tilde{E}^{-n}({\mathcal L})) \subset J_n$.
Let $\lambda_n:u_n(I_n) \to (0,1)$ be the restriction to the
interval $u_n(I_n)$ of an affine map   sending the interval
$u_n(I_n)$ onto the interval $(0,1)$  (see Figure \ref{car09}).
Let   $e_n:(0,1) \to \reals$ be the $C^r$   map defined by $e_n = u
\circ E^n \circ u_n^{-1} \circ \lambda_n^{-1}$.
The  map $e_n$ is the composition of a contraction
$\lambda_n^{-1}$ followed by an expansion $u \circ E^n \circ
u_n^{-1}$. Therefore, by the usual blow-down blow-up technique (see
\cite{pinto:markov_map}), the map $e:(0,1) \to \reals$ given by $e =
\lim_{n\to \infty} e_n$ is a $C^r$ homeomorphism.
Hence, the map
$v_{\mathcal L}:I \to \reals$
defined by $e^{-1} \circ u$ is a solenoidal chart and   is $C^r$
compatible with the charts contained in $U'$. \qed

\subsection{Proof of Theorem \ref{abab}}
\label{usus}

Theorem \ref{abab} follows as a corollary of Lemma \ref{lasda9},
below, by taking $d=2$.

\begin{lemma}
\label{lasda9} The following sets are canonically isomorphic:
\begin{rlist}
\item
The set of all $C^{1+ H\ddot older}$ structures
$U$ for
the expanding circle map $E:S \to S$ of degree $d \geq 2$;
\item The set of all (thca)
solenoids
$(\tilde{E},\tilde{S})$;
\item
The set of all H\"older leaf ratio functions $r:T \to \reals^+$;
\item
The set of all H\"older solenoid functions $s:C \to \reals^+$;
\item The set of all sequences $\{r_0,r_1,\ldots\} \in A(d)$;
\item The set of all affine classes of $d$-adic quasiperiodic    tilings of
the real line that are fixed points of the d-amalgamation operator;
\item The set of all affine classes of  $d$-adic quasiperiodic  fixed grids
of
the real line.
\end{rlist}
\end{lemma}

\noindent {\bf Proof:} The proof of this lemma follows from the
following diagram, where the implications are determined by the
lemmas indicated by their numbers:
\[
\begin{array}{ccccccc}
&   &(i)& \stackrel{\ref{rrr}}{\Longleftarrow}&(ii)&
\stackrel{\ref{lassssspq}}{\Longleftrightarrow} & (vi),(vii)\\
&   &\ref{rrr555}{\Uparrow} & \stackrel{\ref{nnnewww}} {\searrow} &
{\Updownarrow} \ref{lasda}&   &  \\
(v)&\stackrel{\ref{addddp}}{\Longleftrightarrow}&(iv)
&\stackrel{\ref{nnneeewww}}{\Longleftarrow}& (iii) &  &
\end{array}
\]
\qed

\subsection{Proof of Theorem \ref{tata}} \label{jrufcp}

In this section, we use Theorem \ref{th13131313} to prove Theorem
\ref{tata} (we will show Theorem \ref{th13131313} in Section
\ref{pppttt}). In fact, we prove a more general version of Theorem
\ref{tata},
which applies  to expanding circle
maps   of degree $d \geq 2$,
using the following  generalization   of the
  ultra-metric
to the set $C$ of all $d$-adic
integers.
Let $a=\sum_{m=0}^\infty a_m d^m  \in C$ and
$b=\sum_{m=0}^\infty b_m d^m \in C$ be such that
$a_{n} \ldots a_0=b_{n} \ldots b_0$ and $a_{n+1} \ne b_{n+1}$.
For $0 \le i \le n$, let
$A_i = \sum_{m=0}^i a_m d^m$ and
$E_i = \sum_{m=0}^i (d-1) d^m$.
We define the {\it ultra-metric} by
$$|{\bf u}|_s(a,b)=\inf_{0 \le i \le n}
\left\{ 1+  \sum_{j=A_i}^{E_i} \prod_{l=A_i}^j s(l)
+  \sum_{j=0}^{A_i-1}  \prod_{l=j}^{A_i-1}  s(l)
\right\} \ .
$$
Let ${\bf p}$ be the fixed point of the solenoid map
$\tilde{E}$ such that $\pi ({\bf p})$ is the
fixed point  of the expanding circle map
chosen in Section \ref{section2}
to generate the Markov partition of $E$.
Let ${\mathcal L}_{\bf p}$ be the local leaf starting on ${\bf p}$
and ending on its image ${\tilde M} ({\bf p})$  by the monodromy map
${\tilde M}$.
Let $z:J=\pi_S({\mathcal L}_{\bf p}) \to
(0,1)$ be the
corresponding solenoidal chart.
The geometric interpretation of the ultra-metric $|{\bf u}|_s$ is given
by the following equality
$$
|{\bf u}|_s (a,b)=\inf_{0 \le i \le n}
\left\{ |z (I_{a_i \ldots a_0})|  \right\} \ .
$$

\noindent {\bf Proof of Theorem \ref{tata}}.
Let $U$ be a $C^{1+H\ddot older}$ structure for the expanding circle
map $E$, and let $V$ be the
set of all  solenoidal charts  determined by the
(thca) solenoid $(\tilde{E},\tilde{S})_U$.
By Lemma \ref{be}, the set $V$ is contained in $U$  and
the degree of smoothness of
the expanding circle map $E$
when measured in terms of a cover $U'$ of $U$ attains its maximum
when $U' \subset V$.
Let ${\mathcal L}$ and
${\mathcal L'}$ be two local leaves and $u:J=\pi_S({\mathcal L}) \to
\reals$ and $v:J'=\pi_S({\mathcal L'}) \to  \reals$ the
corresponding
  solenoidal charts.
If  $J \cap J' \ne \emptyset$, let
$I_{\beta_1\ldots\beta_n} \subset J \cap J'$ be any   interval at
any level $n$ of the interval partition. Let the points ${\bf x} \in
{\mathcal L}$ and ${\bf y} \in {\mathcal L'}$ be   such that
$\pi_S({\bf x}) = \pi_S({\bf y}) \in S$ is the right endpoint
of the interval $I_{\beta_1\ldots\beta_n}$. Let  $a$ be  the point
$\tilde{\omega} (\tilde{E} ^n({\bf x})) \in C$ and
$b$ the point
$\tilde{\omega} (\tilde{E} ^n({\bf y})) \in C$.
By definition of the metric $|{\bf u}|_s$,
there is a   constant $D_0 = D_0({\mathcal L},{\mathcal L'}) \ge 1$, such that
$D_0^{-1} |z(I_{\beta_1\ldots\beta_n})| \le |{\bf u}|_s (a,b) \le D_0 |z(I_{\beta_1\ldots\beta_n})|$
with respect to the solenoidal chart $z:J=\pi_S({\mathcal L}_{\bf p}) \to
(0,1)$ defined above.
 By Lemma \ref{be},  the overlap maps $z \circ
u^{-1}$ and $z \circ v^{-1}$ are
  $C^{1+H\ddot older}$ smooth.
Therefore, by Lemma \ref{qsq23}, there is a   constant
$D_1 = D_1({\mathcal L},{\mathcal L'}) \ge 1$  such that
\begin{equation}
\label{aak}
D_1^{-1} \le \frac{|{\bf u}|_s
(a,b)}{|u(I_{\beta_1\ldots\beta_n})|} \le D_1  ~~~{\rm  and}~~~ D_1^{-1} \le \frac{|{\bf u}|_s (a,b)}{|v(I_{\beta_1\ldots\beta_n})|} \le D_1 \ .
\end{equation}
Let $I_{\beta_1'\ldots\beta_n'}$ and $I_{\beta_1''
\ldots\beta_n''}$ be adjacent intervals at level $n$ of the
interval partition, such that $I_{\beta_1'\ldots\beta_n'}$ is also
adjacent to $I_{\beta_1\ldots\beta_n}$.
By proof of Lemma \ref{rrr555},
\begin{equation}
\label{bbk}
 s(a)= \frac {|u(I_{\beta_1'\ldots\beta_n'}) |}
{|u(I_{\beta_1\ldots\beta_n})|}, ~~~~~ s(a+1)= \frac
{|u(I_{\beta_1''\ldots\beta_n''})|} {|u(I_{\beta_1'\ldots\beta_n'})|},
\end{equation}
 and
\begin{equation}
\label{cck}
 s(b)= \frac {|v(I_{\beta_1'\ldots\beta_n'}) |}
{|v(I_{\beta_1\ldots\beta_n})|}, ~~~~~ s(b+1)= \frac
{|v(I_{\beta_1''\ldots\beta_n''})|} {|v(I_{\beta_1'\ldots\beta_n'})|}.
\end{equation}
 The interval partition of the expanding circle  map $E$
generates a grid $g_u$ in the set $u(J \cap J')$. Therefore,
using (\ref{aak}), (\ref{bbk}), (\ref{cck}) and
Theorem \ref{th13131313}, the equivalences presented in Tables $2$ and $3$
imply that the
overlap maps $h=v \circ u^{-1}:u(J \cap J') \to v(J \cap J')$
satisfy the equivalences presented in Table $1$.
\qed

\section{Smoothness of diffeomorphisms and cross ratio
distortion of grids} \label{section10}

In the following subsections, we introduce the definitions  of the
degrees of smoothness of a homeomorphism $h:I \to J$ presented in
Tables $2$ and $3$, and we prove the corresponding equivalences between the
degrees of smoothness of $h$ with the ratio   and cross ratio
distortions of intervals contained in a grid of $I$ as presented in
Tables $2$ and $3$. In Section \ref{pppttt}, we   prove Theorem
\ref{th13131313}.

Let $I_\beta$ and $I_{\beta'}$ be two intervals contained in the
real line. We define   the {\it ratio} $r(I_\beta,I_{\beta'})$
between the intervals $I_\beta$  and $I_{\beta'}$   by
$$ r(I_\beta,I_{\beta'})
=\frac{|I_{\beta'}| }{ |I_\beta|} \ .
$$
Let $I_\beta$, $I_{\beta'}$ and $I_{\beta''}$ be contained in the real line,
such that  $I_\beta$ is adjacent to  $I_{\beta'}$,
  and   $I_{\beta'}$ is adjacent to  $I_{\beta''}$.
Recall that the {\it cross ratio} $cr(I_\beta, I_{\beta'}, I_{\beta''})$
is given  by
\begin{eqnarray*}
cr(I_\beta, I_{\beta'}, I_{\beta''})
& = &
\log
\left(
1+
\frac
{|I_{\beta'}|}
{|I_\beta|}
\frac
{|I_\beta|+|I_{\beta'}|+|I_{\beta''}|}
{|I_{\beta''}|}
\right) \ .
\end{eqnarray*}
We note that
$$ cr(I_\beta, I_{\beta'}, I_{\beta''}) =
\log
\left(
(1+r(I_\beta, I_{\beta'}))
(1+r(I_{\beta''}, I_{\beta'}))
\right) \ .
$$
%We note that the cross ratio $cr(I_\beta, I_{\beta'}, I_{\beta''})$
%is   the {\it  Poincare length} of $I_{\beta'}$ in
%$I_\beta \cup I_{\beta'} \cup I_{\beta''}$.

Let $h:I \subset \reals \to J \subset \reals$ be a  homeomorphism, and let
${\mathcal G}_\Omega$
be a grid of the compact interval $I$. We will use the following
  definitions and notations throughout all this section:
\begin{rlist}
\item We will denote
by $J_\beta^n$ the   interval $h(I_\beta^n)$ where $I_\beta^n$ is a grid
interval.
We will denote by $r(n,\beta)$ the ratio $r(I_\beta^n, I_{\beta+1}^n)$
between the grid intervals $I^n_\beta$ and $I^n_{\beta+1}$, and we will
denote
by
$r_h(n,\beta)$ the ratio $r(J_\beta^n, J_{\beta+1}^n)$.
\item
Let $I_\beta$ be an interval contained in $I$  (not necessarily a grid
interval).
The {\it average derivative} $dh(I_\beta)$ is given by
$$
dh(I_\beta)=\frac{|h(I_\beta)|}{|I_\beta|} \ .
$$
We will denote by $dh(n,\beta)$   the average derivative
$dh(I_\beta^n)$  of the grid interval $I_\beta^n$.
\item The {\it logarithmic average derivative} $ldh(I_\beta)$ is given
by
$$
ldh(I_\beta)=\log ( dh (I_\beta)) \ .
$$
We will denote by   $ldh(n,\beta)$  the
logarithmic average derivative $ldh(I_\beta^n)$ of the grid interval
$I_\beta^n$.
\item
Let $I_\beta$ and $I_{\beta'}$  be intervals contained in $I$ (not
necessarily  grid intervals). We recall that the   logarithmic ratio
distortion $lrd (I_\beta, I_{\beta'})$   is given  by $$ lrd
(I_\beta, I_{\beta'})   = \log \left( \frac{|I_{\beta}|}
{|I_{\beta'}|} \frac{|h(I_{\beta'})|}{|h(I_{\beta})|} \right) \ .
$$ Hence, we have $$ lrd (I_\beta, I_{\beta'})
=
\log
  \frac{r(J_\beta, J_{\beta'})}{r(I_\beta, I_{\beta'})}
=
\log \frac{dh(I_{\beta'})}{dh(I_\beta)} \ .
$$
We will denote by $lrd(n,\beta)$ the logarithmic ratio distortion
$lrd(I_\beta^n, I_{\beta+1}^n)$
of the grid intervals  $I_\beta^n$ and $I_{\beta+1}^n$.
\item
Let  the intervals $I_\beta$, $I_{\beta'}$ and $I_{\beta''}$ in $I$
(not necessarily  grid intervals) be such that
$I_\beta$ is adjacent to   $I_{\beta'}$   and
  $I_{\beta'}$ is adjacent to  $I_{\beta''}$.
We recall that the cross ratio distortion  $crd(I_\beta, I_{\beta'},
I_{\beta''})$ is given  by $$ crd (I_\beta, I_{\beta'},
I_{\beta''})= cr(h(I_\beta), h(I_{\beta'}),
h(I_{\beta''}))-cr(I_\beta, I_{\beta'}, I_{\beta''}) \ . $$ We note
that
\begin{equation}
\label{gedrgcdxdd}
crd (I_\beta, I_{\beta'}, I_{\beta''})   =
\log \left( \frac{1+r(h(I_\beta), h(I_{\beta'})) }{1+r(I_\beta,
I_{\beta'})} \frac{ 1+r(h(I_{\beta''}), h(I_{\beta'}))} {
1+r(I_{\beta''}, I_{\beta'})} \right)  \ .
\end{equation}
For all grid intervals $I_\beta^n$,  $I_{\beta+1}^n$ and
$I_{\beta+2}^n$, we will denote by $cr(n,\beta)$ and
$cr_h(n,\beta)$ the cross ratios $cr(I_\beta^n,
I_{\beta+1}^n,I_{\beta+2}^n)$ and
$cr(J_\beta^n, J_{\beta+1}^n,J_{\beta+2}^n)$ respectively.
We will denote by  $crd(n,\beta)$ the cross ratio distortion  given by
$cr_h(n,\beta)-cr(n,\beta)$.
\end{rlist}

\begin{remark}
\label{gridremark} \begin{rlist}
\item[(a)]
We will call properties (vi) and (vii) of a (B,M) grid  ${\mathcal
G}_\Omega$ of an interval $I$, the {\it bounded geometry property}
of the grid.
\item[(b)]
By the bounded geometry property of a (B,M) grid ${\mathcal G}_\Omega$,
there
are constants $0<B_1 <B_2<1$, just depending upon $B$ and $M$, such
that
$$ B_1< \frac{|I_\beta^{n+1}|}{ |I_\alpha^n|} < B_2 \ , $$ for all
$n \ge 1$ and for all grid intervals  $I_\alpha^n$ and
$I_\beta^{n+1}$ such that $I_\beta^{n+1} \subset I_\alpha^n$.
\item[(c)]
We call a $(1,2)$ grid ${\mathcal G}_\Omega$ of $I$  a {\it
symmetric grid of $I$}, i.e. (i) all the intervals at the same level
$n$ have the same length, and (ii) each grid interval at level $n$
is equal to the union of two grid intervals at level $n+1$.
\end{rlist}
\end{remark}

\subsection{Quasisymmetric homeomorphisms}
\label{section11}

The definition of a quasisymmetric homeomorphism that we introduce
in this paper is more     adapted to our problem and, apparently, is
stronger than the usual one,
where   the constant $d$ of  the quasisymmetric condition in Definiton
\ref{quasssi}, below,
is taken to be equal to $1$.
However, in Lemma \ref{rfgerggfgddW1}, we will prove that they are
equivalent.

\begin{defn}
\label{quasssi} Let $d \ge 1$ and  $k  \ge 1$. The  homeomorphism
$h:I \to  J$ satisfies the  {\it $(d,k)$ quasisymmetric condition}
if
\begin{equation}
\label{fgggrrs}
\left|
\log
\frac
{h (x+\delta_2) - h (x)}
{h (x) - h (x-\delta_1)}
\frac
{\delta_1}
{ \delta_2}
\right|
\le \log(k) \ ,
\end{equation}
for all $x-\delta_1,x,x+\delta_2 \in I$ with $\delta_1>0$,
$\delta_2>0$ and $d^{-1} \le \delta_2 \ / \delta_1 \le d$. The
homeomorphism $h$ is {\it   quasisymmetric} if  for every $d \ge 1$
there exists $k_d \ge 1$ such that $h$ satisfies the $(d,k_d)$
quasisymmetric condition.
\end{defn}

\begin{lemma}
\label{lll0}
Let   $h:I \to  J$
be a  homeomorphism and let ${\mathcal G}_\Omega$ be a grid of a compact
interval $I$.
The following statements are equivalent:
\begin{rlist}
\item
The homeomorphism $h:I \to  J$ is
   quasisymmetric.
\item
There is $k({\mathcal G}_\Omega) > 1$ such that
\begin{equation}
\label{sfefwewqesadd} |r_h(n,\beta)| \le k({\mathcal G}_\Omega) \ ,
\end{equation}
for every $n \ge 1$ and every $1 \le \beta \le \Omega(n)$.
\end{rlist}
\end{lemma}

Let  ${\mathcal G}_\Omega$ be a grid of $I$. From Lemma \ref{lll0},
we obtain that a homeomorphism $h:I \to  J$   is   quasisymmetric
if, and only if, the set of all intervals $J_\beta^n$ form a $(B,M)$ grid
for
some $B\ge 1$ and $M > 1$.

\begin{lemma}
\label{rfgerggfgddW1} If, for some $d_0 \ge 1$ and $k_0 \ge 1$, a
homeomorphism $h:I \to  J$ satisfies the $(d_0,k_0)$ quasisymmetric
condition,    then   $h$ is   quasisymmetric.
\end{lemma}

\begin{lemma}
\label{lll1}
Let   $h:I \to  J$
be a  homeomorphism and  ${\mathcal G}_\Omega$ a grid of the compact
interval
$I$.
\begin{rlist}
\item
If $h:I \to  J$ is   quasisymmetric then  there is
  $C_0 \ge  0$ such that
$$
cr_h(n,\beta) \le C_0  \ ,
$$
for every $n \ge 1$ and every $1 \le \beta <\Omega(n)-1$.
\item
If there is
  $C_0 > 1$ such that,  for every $n \ge 1$ and every $1 \le \beta
<\Omega(n)-1$,
$$
cr_h(n,\beta) \le C_0  \ ,
$$
then, for every  closed interval $K$ contained in the interior of
$I$, the homeomorphism $h$ restricted to $K$ is quasisymmetric.
\end{rlist}
\end{lemma}

Before  proving Lemmas \ref{lll0}, \ref{rfgerggfgddW1} and
\ref{lll1}, we will state and prove Lemma \ref{lll33} which we will
use in the proof  of Lemma \ref{lll0} and, later,  in the proof of
Lemma \ref{gfsdgfgd33}.

\begin{lemma}
\label{lll33} Let $\alpha > 1$ and $d \ge 1$. Let   ${\mathcal
G}_\Omega$ be a $(B,M)$ grid of a compact interval $I$. Let
$x-\delta_1, x,x+\delta_2$ contained in $I$ be such that $\delta_1 >
0$, $\delta_2>0$ and  $d^{-1} \le \delta_2/\delta_1 \le d$. Then,
there are intervals $L_1$, $L_2$, $R_1$ and $R_2$ with the following
properties:
\begin{rlist}
\item
\begin{equation}
\label{505050} L_1   \subset    [x-\delta_1,x] \subset L_2~~~{\rm
and}~~~R_1 \subset [x,x+\delta_2] \subset R_2 \ .
\end{equation}
\item
\begin{equation}
\label{tgvcdews}
\alpha^{-1} < \frac{|L_1|}{\delta_1} < \frac{|L_2|}{\delta_1}  <\alpha
~~~{\rm and}~~~
\alpha^{-1} < \frac{|R_1|}{\delta_2} < \frac{|R_2|}{\delta_2}  <\alpha \ .
\end{equation}
\item Let $n_0=n_0(x-\delta_1,x,x+ \delta_2,{\mathcal G}_\Omega)   \ge 1$ be
the
biggest integer such that
$$
[x-\delta_1,x+\delta_2]  \subset I_{\beta}^{n_0} \cup
I_{\beta+1}^{n_0}
$$
for some $1 \le \beta < \Omega(n_0)$. Then, there are integers
$n_1=n_1(\alpha,B,M,d)$
and $n_2=n_2(\alpha,B,M,d)$  such that
\begin{eqnarray*}
L_1 & = &\cup_{i=l+1}^{m-1} I_i^{n_0+n_1}~~~,~~~~L_2=\cup_{i=l}^{m}
I_i^{n_0+n_1} \ , \\
R_1& = &\cup_{i=m+1}^{r-1}  I_i^{n_0+n_1}~~,~~~~R_2=\cup_{i=m}^{r}
I_i^{n_0+n_1} \ ,
\end{eqnarray*}
for some  $l,m,r$ with the property that $l<m<r$ and $r-l \le n_2$.
\end{rlist}
\end{lemma}

\bigskip
\noindent {\bf Proof of Lemma \ref{lll33}.} Let  $0<B_1=B_1(B,M)
<B_2=B_2(B,M) <1$ be  as in Remark \ref{gridremark}.
By construction of $n_0$, there is
$I_\epsilon^{n_0+1}$     with the property  that $I_\epsilon^{n_0+1}
\subset [x-\delta_1,x,x+\delta_2].$ In particular, we have that
either
$I_\epsilon^{n_0+1}\subset  I_{\beta}^{n_0}$
or $I_\epsilon^{n_0+1}\subset  I_{\beta+1}^{n_0}$. Thus, using the
bounded geometry property of a grid  and   Remark \ref{gridremark},
we obtain that
\begin{equation}
\label{fgsggvgfg}
B^{-1}   B_1   |I_{\beta}^{n_0}|  \le    |I_\epsilon^{n_0+1}|    \le
\delta_2
+ \delta_1
   \ .
\end{equation}
Since $d^{-1} \le \delta_2/\delta_1 \le d$, by inequality (\ref{fgsggvgfg}),
we
get
\begin{eqnarray}
\label{grgrg1}
\delta_1  &  \ge  & (1+ D)^{-1} (\delta_2 + \delta_1) \nonumber \\
& \ge &  (1+ D)^{-1} B^{-1}   B_1   |I_{\beta}^{n_0}|   \ .
\end{eqnarray}
Since
$[x-\delta_1,x+\delta_2] \subset I_{\beta}^{n_0} \cup I_{\beta+1}^{n_0}$,
by the  bounded geometry property of a grid, we obtain that
\begin{eqnarray}
\label{grgrg2}
\delta_1 & \le &  |I_{\beta}^{n_0}|+|I_{\beta+1}^{n_0}|  \nonumber  \\
& \le &
(1+B) |I_{\beta}^{n_0}| \ .
\end{eqnarray}
By inequalities (\ref{grgrg1}) and (\ref{grgrg2}), there is $A=A(B_0,B_1,d)
> 1$
such that
\begin{equation}
\label{grgrg3}
  A^{-1}  |I_{\beta}^{n_0}| \le  \delta_1   \le    A |I_{\beta}^{n_0}| \ .
\end{equation}
Similarly, we have
\begin{equation}
\label{grgrg4}
A^{-1}  |I_{\beta}^{n_0}| \le  \delta_2   \le    A |I_{\beta}^{n_0}| \ .
\end{equation}
Take $0<\theta(\alpha) < 1$   such that $\alpha^{-1} \le 1-\theta <1+\theta
\le
\alpha$. Let
$n_1=n_1(B,B_2,A,\theta)$ be the smallest integer such that
\begin{equation}
\label{grgrg45} B_2^{n_1} \le  B^{-1} \theta A^{-1} /2 \ .
\end{equation}
Let  $l < m < r$ be such that
$x-\delta_1 \in  I_{l}^{n_0+n_1}$, $x \in  I_{m}^{n_0+n_1}$
and $x+\delta_2 \in  I_{r}^{n_0+n_1}$. Then, by the bounded geometry
property of a grid, there is $n_2=2 M n_1 \ge 1$ such that $r-l \le
n_2$. Hence, the intervals
\begin{eqnarray*}
L_1 & = &\cup_{i=l+1}^{m-1} I_i^{n_0+n_1}~~~,~~~~L_2=\cup_{i=l}^{m}
I_i^{n_0+n_1} \ , \\
R_1& = &\cup_{i=m+1}^{r-1}  I_i^{n_0+n_1}~~,~~~~R_2=\cup_{i=m}^{r}
I_i^{n_0+n_1}
\end{eqnarray*}
satisfy property (i) and property (iii) of Lemma \ref{lll33}. Let us
prove that the intervals $L_1$, $L_2$, $R_1$ and $R_2$   satisfy
property (ii) of Lemma \ref{lll33}. By the bounded geometry property
of a grid and inequality \ref{grgrg45}, we get
\begin{equation}
\label{fgghbnhg}
|I_i^{n_0+n_1}| \le    B  B_2^{n_1} |I_{\beta}^{n_0}| \le \theta A^{-1}
|I_{\beta}^{n_0}|/ 2 \ ,
\end{equation}
for all $l \le i \le r$.
Thus, by inequalities (\ref{grgrg3}) and (\ref{fgghbnhg}), we get
\begin{eqnarray}
\label{defefef1}
|L_1|/\delta_1 & \ge & (\delta_1 - |I_l^{n_0+n_1}| -|I_m^{n_0+n_1}|)/
\delta_1
\nonumber \\
& \ge & (\delta_1 -      \theta A^{-1} |I_{\beta}^{n_0}|) / \delta_1
\nonumber
\\
& \ge &   1-  \theta   \ .
\end{eqnarray}
Again, by inequalities (\ref{grgrg3}) and (\ref{fgghbnhg}), we get
\begin{eqnarray}
\label{defefef333}
|L_2|/\delta_1 & \le & (\delta_1 + |I_l^{n_0+n_1}| +|I_m^{n_0+n_1}|)/
\delta_1
\nonumber \\
& \le & (\delta_1 +      \theta
A^{-1} |I_{\beta}^{n_0}|) / \delta_1 \nonumber \\
& \le &   1+   \theta  \ .
\end{eqnarray}
Similarly, using  inequalities  (\ref{grgrg4}) and (\ref{fgghbnhg}), we
obtain
that
\begin{equation}
\label{defefef2}
|R_1|/ \delta_2   \ge    1-  \theta
~~~{\rm and}~~~
|R_2|/ \delta_2    \le     1+   \theta   \ .
\end{equation}
Noting that $\alpha^{-1} \le 1-\theta <1+\theta \le \alpha$ and putting
together
inequalities
(\ref{defefef1}), (\ref{defefef333}) and     (\ref{defefef2}), we obtain
that
the intervals
$L_1$, $L_2$, $R_1$ and
$R_2$   satisfy property (ii) of Lemma \ref{lll33}.
\qed

\bigskip
\noindent {\bf Proof of Lemma \ref{lll0}.} {\it Let us prove that
statement (i) implies statement (ii)}. For every level $n \ge 1$ and
every $1 \le \beta < \Omega(n)$, let $x-\delta_1,x,x+\delta_2 \in I$
be such that $I_\beta^n=[x-\delta_1,x]$ and
$I_{\beta+1}^n=[x,x+\delta_2]$. Hence,
$$
\frac{r_h(n,\beta)}{r(n,\beta)}=
\frac
{h (x+\delta_2) - h (x)} {h (x) - h (x-\delta_1)} \ .
$$
Since  $h:I \to  J$ is
  $(k,B)$ quasisymmetric, for some $k=k(B)$, we have
$$
k^{-1} <
\frac
{h (x+\delta_2) - h (x)} {h (x) - h (x-\delta_1)}
\frac
{\delta_1}{ \delta_2} <  k\ ,
$$
and so, we get
\begin{equation}
\label{Wgfbgrgrg}
k^{-1} <  r_h(n,\beta)/r(n,\beta)  <  k \ .
\end{equation}
Since, by the bounded geometry property of a grid  ${\mathcal
G}_\Omega$, there is $B \ge 1$ such that $ B^{-1} \le r(n,\beta) \le
B $,
we get $k^{-1}B^{-1} \le r_h(n,\beta) \le kB$.

\bigskip

\noindent {\it Let us prove that statement (ii) implies statement
(i)}. Let $B \ge  1$ and $M>1$ be as in the bounded geometry
property of a grid. Let $d \ge 1$. Let  $x-\delta_1, x,
x+\delta_2\in I$, be such that $\delta_1 > 0$, $\delta_2 >0$ and
$d^{-1} \le \delta_2 / \delta_1 \le d$. Let $L_1$, $L_2$, $R_1$ and
$R_2$ be the intervals as constructed in Lemma \ref{lll33} with the
constant $\alpha=2$ in Lemma \ref{lll33}. Hence, we have that
\begin{eqnarray*}
\label{tgvcdews1}
|L_1| & = &  |I_l^{n_0+n_1}| \left ( 1+ \sum_{i=l+1}^{m-2} \prod_{j=l}^{i}
r(n_0+n_1,j) \right ) \ , \\
|L_2| & = & |I_l^{n_0+n_1}|  \left (1+ \sum_{i=l}^{m-1} \prod_{j=l}^{i}
r(n_0+n_1,j) \right) \ , \\
|R_1| & = &  |I_l^{n_0+n_1}|      \left( \sum_{i=m}^{r-2} \prod_{j=l}^{i}
r(n_0+n_1,j) \right) \ , \\
|R_2| & = & |I_l^{n_0+n_1}|   \left( \sum_{i=m-1}^{r-1} \prod_{j=l}^{i}
r(n_0+n_1,j) \right) \ .
\end{eqnarray*}
Hence,
by monotonicity of the homeomorphism $h$, we obtain that
\begin{equation}
\label{fgrergb0}
\frac{|h(R_1)|}{|h(L_2)|}
\frac{|L_1|}{|R_2|}
\le
\frac
{h (x+\delta_2) - h (x)} {h (x) - h (x-\delta_1)}
\frac
{\delta_1}{ \delta_2}
\le
\frac{|h(R_2)|}{|h(L_1)|}
\frac{|L_2|}{|R_1|} \ .
\end{equation}
Since, by the bounded geometry property of a grid,    $B^{-1} <
r(n_0+n_1,j) < B$ and, by Lemma \ref{lll33}, $l < m <r$ and $r-l \le
n_2(B,M,d)$, we get that there is $C_1=C_1(B,n_2) > 1$ such that
\begin{eqnarray}
\label{ggrgrgggbv10}
C_1^{-1} & \le & \frac{|L_1|}{|R_2|} =
\frac
{ 1+ \sum_{i=l+1}^{m-2} \prod_{j=l+1}^{i} r(n_0+n_1,j)}
{\sum_{i=m-1}^{r-1} \prod_{j=l+1}^{i} r(n_0+n_1,j) }
\le C_1 \ ,  \nonumber  \\
C_1^{-1} & \le & \frac{|L_2|}{|R_1|}  =
\frac
{1+ \sum_{i=l}^{m-1} \prod_{j=l}^{i} r(n_0+n_1,j) }
{ \sum_{i=m}^{r-2} \prod_{j=l}^{i} r(n_0+n_1,j) }
\le C_1 \ .
\end{eqnarray}
By inequality (\ref{sfefwewqesadd}) of statement (ii), there is
$k=k({\mathcal G}_\Omega)
>1$ such that $k^{-1} < r_h(n_0+n_1,j) < k$ for every $1\le j <
\Omega(n_0+n_1)$. Hence,   there is $C_2=C_2(k,n_2) >1$ such that
\begin{eqnarray}
\label{hbgwffeeeb0}
C_2^{-1} &  \le &  \frac{|h(R_1)|}{|h(L_2)|}
=
\frac
{\sum_{i=m}^{r-2} \prod_{j=l}^{i} r_h(n_0+n_1,j)}
{1+ \sum_{i=l}^{m-1} \prod_{j=l}^{i} r_h(n_0+n_1,j)} \le C_2 \ ,   \nonumber
\\
C_2^{-1} &  \le & \frac{|h(R_2)|}{|h(L_1)| } =
\frac
{\sum_{i=m}^{r-2} \prod_{j=l}^{i} r_h(n_0+n_1,j)}
{1+ \sum_{i=l}^{m-1} \prod_{j=l}^{i} r_h(n_0+n_1,j)}
\le C_2 \ .
\end{eqnarray}
Putting together equations (\ref{fgrergb0}),  (\ref{ggrgrgggbv10})
and   (\ref{hbgwffeeeb0}), we obtain that
$$
C_1^{-1} C_2^{-1}
\le
\frac{|h(R_1)|}{|h(L_2)|}
\frac{|L_1|}{|R_2|}
\le
\frac
{h (x+\delta_2) - h (x)} {h (x) - h (x-\delta_1)}
\frac
{\delta_1}{ \delta_2}
\le
\frac{|h(R_2)|}{|h(L_1)|}
\frac{|L_2|}{|R_1|}
\le C_1C_2 \ .
$$
\qed

\bigskip
\noindent {\bf Proof of Lemma \ref{rfgerggfgddW1}:} If a
homeomorphism  $h:I \to J$ satisfies   the $(d_0,k_0)$
quasisymmetric condition for some $d_0 \ge 1$ and $k_0 \ge 1$ then
$h$ satisfies statement (ii) of Lemma \ref{lll0} with respect to a
symmetric grid (see definition of a symmetric grid in Remark
\ref{gridremark}). Hence, by statement (i) of Lemma \ref{lll0}, the
homeomorphism  $h$ is quasisymmetric. \qed

\bigskip
\noindent {\bf Proof of Lemma \ref{lll1}.} {\it Let us prove
statement  (i)}. By Lemma \ref{lll0},
  there is $C_1 \ge 1$ such that
$C^{-1}_1 \le r_h(n,\beta) \le C_1$
for every level $n$ and every
$1 \le \beta < \Omega(n)$.
Therefore, there is $C_2 >  0$
such that,   for every level $n$ and every
$1 \le \beta < \Omega(n)-1$,
\begin{equation}
\label{grggvs} |cr_h(n,\beta)| = \left| \log \left (
(1+r_h(n,\beta)) (1+r_h(n,\beta+1))^{-1} \right) \right | \le C_2  \
.
\end{equation}

\bigskip

\noindent {\it Let us prove statement (ii)}. By the  bounded
geometry property of a grid, there is  $n_0 \ge 1$ large enough such
that the grid intervals $I_1^{n_0}$ and $I_{\Omega(n)-1}^{n_0}$ do
not intersect the interval $L$. The grid ${\mathcal G}_\Omega$ of
$I$ induces, by restriction, a grid of the interval $L'=
\cup_{\beta=2}^{\Omega(n)-2} I_\beta^{n_0}$ which contains $L$.
Hence, by  Lemma \ref{lll0}, it is enough to prove that there is
$C_1 \ge 1$ such that $C_1^{-1} \le r_h(n,\beta) \le C_1$ for  every
grid interval $I_\beta^n \subset L'$. Now, we will consider
separately  the following two possible cases: either (i)
$r_h(n,\beta) \le 1$ or (ii) $r_h(n,\beta) > 1$.

\bigskip
\noindent
{\it Case (i)}. Let   $r_h(n,\beta)=|J_{\beta+1}^n|/|J_\beta^n|\le 1$.
By hypotheses of statement (ii), there is
$C_2 > 1$ such that
$$
cr_h(n,\beta-1) =
\log \left(1+
\frac
{|J_{\beta}^n|}
{|J_{\beta-1}^n|}
\frac
{|J_{\beta-1}^n|+|J_{\beta}^n|+|J_{\beta+1}^n|}
{|J_{\beta+1}^n|}
\right)
\le C_2 \ .
$$
Hence, there is
$C_3 > 1$ such that
$$
1 \le
\frac
{|J_\beta^n|}
{|J_{\beta+1}^n|}
\le
\frac
{|J_{\beta}^n|}
{|J_{\beta+1}^n|}
\frac
{|J_{\beta-1}^n|+|J_{\beta}^n|+|J_{\beta+1}^n|}
{|J_{\beta-1}^n|}
\le C_3 \ ,
$$
and so $C_3^{-1}  \le r_h(n,\beta) \le 1$.

\bigskip
\noindent
{\it Case (ii)}.  Let   $r_h(n,\beta)= |J_{\beta+1}^n|/|J_\beta^n| > 1$. By
hypotheses, there is
$C_2 > 1$ such that
$$
cr_h(n,\beta) =
\log \left(1+
\frac
{|J_{\beta+1}^n|}
{|J_\beta^n|}
\frac
{|J_\beta^n|+|J_{\beta+1}^n|+|J_{\beta+2}^n|}
{|J_{\beta+2}^n|}
\right)
\le C_2 \ .
$$
Hence, there is
$C_3 > 1$ such that
$$
1 \le
\frac
{|J_{\beta+1}^n|}
{|J_\beta^n|}
\le
\frac
{|J_{\beta+1}^n|}
{|J_\beta^n|}
\frac
{|J_\beta^n|+|I_{\beta+1}^n|+|I_{\beta+2}^n|}
{|J_{\beta+2}^n|}
\le C_3 \ ,
$$
and so $1 <r_h(n,\beta) \le C_3$.
\qed

\subsection{Horizontal and vertical translations of ratio distortions}

  Lemmas \ref{equalities} and \ref{equalities3}   are the key to understand
the
relations
between ratio and cross ratio distortions. We will use them in the
following subsections.

\begin{lemma}
\label{equalities} Let $h:I \subset \reals \to J \subset \reals$ be
a quasisymmetric homeomorphism and ${\mathcal G}_\Omega$ a grid of
the closed interval $I$. Then, the logarithmic ratio distortion and
cross ratio distortion satisfy the following estimates:
\begin{eqnarray}
\label{batatas2}
\frac{r_h(n,\beta)}{r(n,\beta)} & \in &
1+ lrd(n,\beta)   \pm \cO (lrd(n,\beta)^2) \\
\label{batatas222}
\frac{r(n,\beta)}{r_h(n,\beta)} & \in &
1- lrd(n,\beta)   \pm \cO (lrd(n,\beta)^2)
\end{eqnarray}
\begin{eqnarray}
crd(n,\beta) & \in &
\frac{lrd(n,\beta)}{1+r(n,\beta)^{-1}}
-
\frac{lrd(n,\beta+1)}{1+r(n,\beta+1)}
\pm \cO (lrd(n,\beta)^2,lrd(n,\beta+1)^2)  \nonumber \\
\label{batatas3}
&=&
\frac{|I_{\beta+1}^n|lrd(n,\beta)}{|I_{\beta}^n|+|I_{\beta+1}^n|}
-
\frac{|I_{\beta+1}^n|lrd(n,\beta+1)}{|I_{\beta+1}^n|+|I_{\beta+2}^n|}
\pm \cO (lrd(n,\beta)^2,lrd(n,\beta+1)^2) \ .
\end{eqnarray}
\end{lemma}

In what follows, we will use the following notations:
\begin{eqnarray*}
\label{dfffedede2}
L_1(n,\beta,p) & = & \max_{0 \le i\le p} \{lrd(n,\beta+i)^2 \}  \nonumber \\
L_2(n,\beta,p)  & = &  \max_{0 \le i_1 \le i_2 < p} \{|lrd(n,\beta+i_1)
lrd(n,\beta+i_2)| \} \nonumber \\
C(n,\beta,p) & = &  \max_{0 \le i  < p} \{|crd(n,\beta+i)| \} \nonumber  \\
M_1(n,\beta,p) & = & \max \{L_1(n,\beta,p), C(n,\beta,p) \} \nonumber \\
M_2(n,\beta,p)  & = &  \max \{L_2(n,\beta,p), C(n,\beta,p) \} \ . \nonumber
\end{eqnarray*}

\begin{lemma}
\label{equalities3} Let $h:I \subset \reals \to J \subset \reals$ be
a  quasisymmetric homeomorphism and let ${\mathcal G}_\Omega$ be a
grid of the closed interval $I$. Then, the logarithmic ratio
distortion and the cross ratio distortion satisfy the following
estimates:
\begin{rlist}
\item ({\it lrd-horizontal translations}) There is a constant $C(i)>0$,
not
depending upon the level $n$ and
  not depending  upon $1 \le \beta \le \Omega(n)$, such that
\begin{eqnarray}
\label{batatas4}
lrd(n,\beta+i) & \in &
\left( \prod_{k=0}^{i-1}
r(n,\beta+k) \right )
\frac
{ 1+r(n,\beta+i)}
{1+r(n,\beta)}
lrd(n,\beta)
\pm C(i) M_1(n,\beta,i) \nonumber \\
&= & \frac{|I_{\beta+i}^n|+|I_{\beta+i+1}^n|}{|I_{\beta}^n|+|I_{\beta+1}^n|}
lrd(n,\beta) \pm   C(i) M_1(n,\beta,i) \ .
\end{eqnarray}
\item ({\it lrd-vertical translations})
Let $I_{\alpha}^{n-1}$ and $I_{\alpha+1}^{n-1}$ be two adjacent grid
intervals.
Take $\beta=\beta(n,\alpha)$ and $p=p(n,\alpha)$
such that  $I_{\beta}^n,\ldots,I_{\beta+p}^n$ are all the grid intervals
contained in the union
$I_{\alpha}^{n-1} \cup I_{\alpha+1}^{n-1}$.
Then, for every $0 \le i< p$, we have
\begin{equation}
\label{batatas5}
lrd(n-1,\alpha)   \in
\frac{|I_{\alpha}^{n-1}|+|I_{\alpha+1}^{n-1}|}{|I_{\beta+i}^n|+|I_{\beta+i+1}^n|}
lrd(n,\beta+i) \pm \cO( M_2(n,\beta,p) ) \ .
\end{equation}
\end{rlist}
\end{lemma}

\bigskip
\noindent {\bf Proof of Lemma \ref{equalities}:} \noindent {\it Let
us prove inequality (\ref{batatas2}).} By    Taylor series
expansion, we have that $\log(x)=x-1  \pm \cO((\log x)^2)$ for every
$x$ in a small neighbourhood of $1$. Hence, using that $h$ is
quasisymmetric, we get
\begin{eqnarray*}
lrd(n,\beta) & = & \log \left( \frac{r_h(n,\beta)}{r(n,\beta)}
\right) \\ & \in &
  \frac{r_h(n,\beta)}{r(n,\beta)} -1     \pm \cO(lrd(n,\beta)^2)  \ .
\end{eqnarray*}

\bigskip
\noindent {\it Let us prove inequality (\ref{batatas222}).} By
Taylor series expansion, we have that $1/(1+x) \in 1-x  \pm
\cO(x^2)$ for every $x$ in a small neighbourhood of $0$. Thus, using
that $h$ is quasisymmetric, we obtain that
\begin{eqnarray*}
\label{ghbrtbrc1}
\frac{r(n,\beta)}{r_h(n,\beta)} & = &
\frac{1}{1+ r_h(n,\beta) r(n,\beta)^{-1} -1} \nonumber \\
& \in &
1- \left( \frac{r_h(n,\beta)}{r(n,\beta)} -1  \right)
\pm \cO \left(\left( \frac{r_h(n,\beta)}{r(n,\beta)} -1  \right)^2 \right) \
.
\end{eqnarray*}
Hence, using inequality (\ref{batatas2}), we get
$$
\frac{r(n,\beta)}{r_h(n,\beta)}  \in
1-
lrd(n,\beta) \pm \cO(lrd(n,\beta)^2) \ .
$$

\bigskip
\noindent {\it Let us prove inequality (\ref{batatas3}).} By
definition of cross ratio distortion, we have
$$
crd(n,\beta) =
\log
\frac{1+r_h(n,\beta)}
{1+r(n,\beta)}
+ \log \frac{1+r_h(n,\beta+1)^{-1}}
{1+r(n,\beta+1)^{-1}} \ .
$$
By    Taylor series expansion, we have that $\log(x+1)=x   \pm
\cO(x^2)$ for every $x$ in a small neighbourhood of $0$. By the
bounded geometry property of a grid, there is $C>1$ such that
$C^{-1} \le  1+r(n,\beta)^{-1} \le C$ for every level $n\ge1$ and
$\beta=1,\ldots,\Omega(n)$. Hence,  using inequality
(\ref{batatas2}), we get
\begin{eqnarray}
\label{ovos111}
\log
\frac{1+r_h(n,\beta)}
{1+r(n,\beta)}
&  = &
\log \left(1+
\frac{r_h(n,\beta) r(n,\beta)^{-1} -1}
{1+r(n,\beta)^{-1}} \right) \nonumber \\
& \in &
\frac{lrd(n,\beta))}
{1+r(n,\beta)^{-1}} \pm \cO(lrd(n,\beta)^2)   \ .
\end{eqnarray}
Similarly,  using inequality (\ref{batatas222}), we obtain
\begin{eqnarray}
\label{ovos333}
\log
\frac{1+r_h(n,\beta+1)^{-1}}
{1+r(n,\beta+1)^{-1}}
&  = &
\log \left(1+
\frac{r(n,\beta+1) r_h(n,\beta+1)^{-1} -1}
{1+r(n,\beta+1)} \right) \nonumber \\
& \in &
\frac{-lrd(n,\beta+1)}
{1+r(n,\beta+1)}
\pm
\cO(lrd(n, \beta+1)^2) \ .
\end{eqnarray}
Putting together equations (\ref{ovos111}) and  (\ref{ovos333}), we get
\begin{eqnarray*}
crd(n,\beta) & = &
\log
\frac{1+r_h(n,\beta)}
{1+r(n,\beta)}
+
\log
\frac{1+r_h(n,\beta+1)^{-1}}
{1+r(n,\beta+1)^{-1}} \\
& \in &
\frac{lrd(n,\beta))}
{1+r(n,\beta)^{-1}}
-
\frac{lrd(n,\beta+1)}
{1+r(n,\beta+1)}
\pm
\cO(lrd(n, \beta)^2,lrd(n, \beta+1)^2) \  .
\end{eqnarray*}
\qed

\bigskip
\noindent
{\bf Proof of Lemma \ref{equalities3}.}
{\it Let us prove inequality (\ref{batatas4}).}
Using inequality (\ref{batatas3}), we get
\begin{eqnarray*}
lrd(n,\beta+i+1)
& \in &
lrd(n,\beta+i)  \frac
{ 1+r(n,\beta+i+1)}
{1+r(n,\beta+i)^{-1}}
\pm \cO (M_1(n,\beta+i,1) )  \\
& \subset &
lrd(n,\beta+i) r(n,\beta+i) \frac
{ 1+r(n,\beta+i+1)}
{1+r(n,\beta+i)}
\pm \cO (M_1(n,\beta+i,1) ) \ .
\end{eqnarray*}
Hence, we obtain
\begin{eqnarray*}
lrd(n,\beta+i) & \in &
lrd(n,\beta)
\prod_{k=0}^{i-1}
\left (
  r(n,\beta+k) \frac
{ 1+r(n,\beta+k+1)}
{1+r(n,\beta+k)}
\right)
\pm C(i)  M_1(n,\beta,i)  \\
& \subset &
lrd(n,\beta)
\frac
{ 1+r(n,\beta+i)}
{1+r(n,\beta)}
\prod_{k=0}^{i-1}
r(n,\beta+k)
\pm C(i) M_1(n,\beta,i) \ ,
\end{eqnarray*}
where the constant $C(i) > 0$ does not depend upon $n$ and upon $1 \le \beta
\le
\Omega(n)$.
Since
$$
\frac
{1+r(n,\beta+i)}
{1+r(n,\beta)}
\prod_{k=0}^{i-1} r(n,\beta+k)  =
\frac{|I_{\beta+i}^n|+|I_{\beta+i+1}^n|}
{|I_{\beta}^n|+|I_{\beta+1}^n|} \ ,
$$
  we get
\begin{eqnarray*}
lrd(n,\beta+i) & \in &
lrd(n,\beta)
\frac
{ 1+r(n,\beta+i)}
{1+r(n,\beta)}
\prod_{k=0}^{i-1}
r(n,\beta+k)
\pm C(i) M_1(n,\beta,i)  \\
& = & lrd(n,\beta)
\frac{|I_{\beta+i}^n|+|I_{\beta+i+1}^n|}
{|I_{\beta}^n|+|I_{\beta+1}^n|}
\pm C(i) M_1(n,\beta,i)  \ .
\end{eqnarray*}

\bigskip
\noindent {\it Let us prove inequality (\ref{batatas5}).}
Let  $0 <m=  m(n,\alpha) < p$ be such that
$I_{\beta}^n,\ldots,I_{\beta+m}^n$ are all the grid intervals contained
in
$I_{\alpha}^{n-1}$  and
$I_{\beta+m+1}^n,\ldots,I_{\beta+p}^n$ are all the grid intervals contained
in
$I_{\alpha}^{n-1}$.
For simplicity of exposition, we introduce the following definitions:
\begin{rlist}
\item We define $a_0=0$, $a_{h,0}=0$ and,
for every $0 < j <p$, we define
$$
a_j= \frac{|I_{\beta+j}^n|}{|I_{\beta}^n|}=\prod_{i=0}^{j-1} r(n,\beta+i)
~~~{\rm and}~~~
a_{h,j}= \frac{|J_{\beta+j}^n|}{|J_{\beta}^n|}=\prod_{i=0}^{j-1}
r_h(n,\beta+i)
\ .
$$
\item We define
$$
R  =   \frac{ |I_{\alpha}^{n-1}|}{|I_{\beta}^n|} ~,~~~~
R' =  \frac{|I_{\alpha+1}^{n-1}|}{|I_{\beta}^n|} ~,~~~~
R_h  =  \frac{|J_{\alpha}^{n-1}|}{|J_{\beta}^n|} ~,~~~~
R_h'  = \frac{|J_{\alpha+1}^{n-1}|}{|J_{\beta}^n|} \ .
$$
Thus,
$$
R   =  \sum_{j=0}^{m-1} a_j ~,~~~~
R'  =  \sum_{j=m}^{p-1} a_j ~,~~~~
R_h =  \sum_{j=0}^{m-1} a_{h,j} ~,~~~~~
R_h'=\sum_{j=m}^{p-1} a_{h,j} \ .
$$
\item We define
$$ E   =
\sum_{j=1}^{m-1}
a_j
\left(
\sum_{i=0}^{j-1} lrd(n,\beta+i)
\right)   ~~~{\rm and}~~~
E'  =
\sum_{j=m}^{p-1}
a_j
\left(
\sum_{i=0}^{j-1}   lrd(n,\beta+i)
\right)
\ .
$$
\end{rlist}
We will separate the proof of inequality (\ref{batatas5})  in three parts.
In the first part, we will prove that
\begin{equation}
\label{lalap}
lrd(n-1,\alpha)   \in
\frac{E'}{R'}
-
\frac{E}{R}  \pm \cO(L_2(n,\beta,p)) \ .
\end{equation}
In the second part, we will prove that
\begin{equation}
\label{lalap2}
\frac{E'}{R'} - \frac{E}{R}
\in
lrd(n,\beta) \frac{|I_{\alpha}^n|+|I_{\alpha+1}^n|}
{|I_{\beta}^n|+|I_{\beta+1}^n|}
\pm   \cO(M_1(n,\beta,p))  \ .
\end{equation}
In the third part, we will use the previous parts to prove  inequality
(\ref{batatas5})  in the case where
$i=0$. Then, we will use inequality (\ref{batatas4}) to
extend, for every $0 \le i < p$, the proof of inequality (\ref{batatas5}).

\bigskip
\noindent {\it First part.}
By inequality (\ref{batatas2}), we have that
$$r_h(n,\beta+i)\in r(n,\beta+i)(1+lrd(n,\beta+i)) \pm
\cO(lrd((n,\beta+i)^2) \ .
$$
Hence, for every $1 \le j < p$, we get
\begin{eqnarray*}
a_{h,j} & = & \prod_{i=0}^{j-1} r_h(n,\beta+i)\\
& \in & \prod_{i=0}^{j-1} \left( r(n,\beta+i)(1+lrd(n,\beta+i)) \pm
\cO(lrd((n,\beta+i)^2) ) \right) \\
& \subset &   \prod_{i=0}^{j-1} r(n,\beta+i)
\left(1+ \sum_{i=0}^{j-1}lrd(n,\beta+i)  \pm
\cO( L_2(n,\beta+1,j))  \right) \\
& \subset &   a_j + a_j \sum_{i=0}^{j-1} lrd(n,\beta+i)   \pm
\cO\left(a_j  L_2(n,\beta+1,j)\right) \ .
\end{eqnarray*}
Thus,
\begin{eqnarray}
\label{gerreg1}
R_h & = & \sum_{j=0}^{m-1} a_{h,j} \nonumber \\
&\in &  \sum_{j=0}^{m-1}
a_j +  \sum_{j=1}^{m-1} a_j \sum_{i=0}^{j-1} lrd(n,\beta+i)   \pm
\cO\left( \sum_{j=0}^{m-1} a_j  L_2(n,\beta,j)\right) \nonumber   \\
&\subset & R +  E \pm  \cO(R  L_2(n,\beta,m)))   \ .
\end{eqnarray}
Similarly, we have
\begin{eqnarray}
\label{gerreg2}
R_h' & = & \sum_{j=m}^{p-1} a_{h,j} \nonumber \\
&\in &  \sum_{j=m}^{p-1}
a_j +  \sum_{j=m}^{p-1} a_j \sum_{i=0}^{j-1} lrd(n,\beta+i)   \pm
\cO \left( \sum_{j=m}^{n-1} a_j  L_2(n,\beta,j) \right) \nonumber  \\
&\subset & R' +  E' \pm  \cO(R'  L_2(n,\beta,p))   \ .
\end{eqnarray}
By  inequalities (\ref{gerreg1}) and (\ref{gerreg2}), we obtain that
\begin{eqnarray*}
lrd(n-1,\alpha) & = &
\log  \frac {R_h'}{R'}
\frac {R}{R_h} \\
& \in  &
\log  \frac{R'+E'\pm\cO(R'  L_2(n,\beta,p))}{R'}
-
\log  \frac{R+E\pm \cO(R  L_2(n,\beta,m)) }{R} \\
& \subset  &
\frac{E'}{R'} - \frac{E}{R} \pm \cO( L_2(n,\beta,p)) \ .
\end{eqnarray*}

\bigskip
\noindent {\it Second part.}
By inequality (\ref{batatas4}), for every $1 \le j <p$, we obtain
\begin{eqnarray*}
\sum_{i=0}^{j-1} lrd(n,\beta+i)  & \in  &
\sum_{i=0}^{j-1}  \left( \frac{a_i(1+r(n,\beta+i))}{1+r(n,\beta)}
lrd(n,\beta)
\pm
\cO(M_1(n,\beta,i)) \right)  \\
& \subset  &
\frac{lrd(n,\beta)}{1+r(n,\beta)} \sum_{i=0}^{j-1} (a_{i} + a_{i+1}) \pm
\cO(M_1(n,\beta,j)) \ .
\end{eqnarray*}
Hence, we obtain that
\begin{eqnarray}
\label{rfgerfda1}
E & = &
\sum_{j=1}^{m-1}
a_j   \sum_{i=0}^{j-1} lrd(n,\beta+i)   \nonumber \\
& \in  &
\sum_{j=1}^{m-1}
a_j  \left(
\frac{lrd(n,\beta)}{1+r(n,\beta)} \sum_{i=0}^{j-1} (a_{i} + a_{i+1}) \pm
\cO(M_1(n,\beta,j))
\right)  \nonumber \\
& \subset  &
\frac{lrd(n,\beta)}{1+r(n,\beta)} \sum_{j=1}^{m-1}
   a_j   \sum_{i=0}^{j-1} (a_{i} + a_{i+1})
\pm  \cO \left(\sum_{j=1}^{m-1} a_j  M_1(n,\beta,j)\right)  \nonumber \\
& \subset &
\frac{lrd(n,\beta)}{1+r(n,\beta)}
R (a_1+\ldots + a_{m-1}) \pm     \cO(R  M_1(n,\beta,m)) \ .
\end{eqnarray}
Similarly, we have
\begin{eqnarray}
\label{rfgerfda2} E' & = & \sum_{j=m}^{p-1} a_j \sum_{i=0}^{j-1}
lrd(n,\beta_i)    \nonumber \\ & \in &
\frac{lrd(n,\beta)}{1+r(n,\beta)} \sum_{j=m}^{p-1} a_j
\sum_{i=0}^{j-1} (a_{i} + a_{i+1}) \pm \cO\left(\sum_{j=m}^{p-1}
a_j M_1(n,\beta,j) \right)  \nonumber \\ & \subset &
\frac{lrd(n,\beta)}{1+r(n,\beta)} R' (1+2a_1+\ldots + 2
a_{m-1}+a_m +\ldots +a_{p-1})  \nonumber \\ && \pm \cO(R'
M_1(n,\beta,p)) \ .
\end{eqnarray}
Putting together   inequalities (\ref{rfgerfda1}) and (\ref{rfgerfda2}), we
obtain that
\begin{eqnarray*}
\frac{E'}{R'} - \frac{E}{R}
& \in &
\frac{lrd(n,\beta)}{1+r(n,\beta)}  (1+a_1+\ldots +a_{p-1})
  \pm   \cO( M_1(n,\beta,p)) \\
& \subset &
lrd(n,\beta)
\frac{|I_{\alpha}^n|+|I_{\alpha+1}^n|}
{|I_{\beta}^n|+|I_{\beta+1}^n|}
\pm   \cO( M_1(n,\beta,p)) \ .
\end{eqnarray*}

\bigskip
\noindent {\it Third part.}
In the case where $i=0$,
inequality (\ref{batatas5}) follows, from
putting together inequalities (\ref{lalap}) and (\ref{lalap2}), since
\begin{eqnarray*}
lrd(n-1,\alpha)   & \in &
\frac{E'}{R'}
-
\frac{E}{R}  \pm \cO(L_2(n,\beta,p))  \\
&\subset & lrd(n,\beta) \frac{|I_{\alpha}^n|+|I_{\alpha+1}^n|}
{|I_{\beta}^n|+|I_{\beta+1}^n|}
\pm   \cO(M_2(n,\beta,p))  \ .
\end{eqnarray*}
By inequality (\ref{batatas4}), for every $0 < i < p$,
we have
$$
\frac{|I_{\alpha}^{n-1}|+|I_{\alpha+1}^{n-1}|}{|I_{\beta}^n|+|I_{\beta+1}^n|}
lrd(n,\beta) \in
\frac{|I_{\alpha}^{n-1}|+|I_{\alpha+1}^{n-1}|}{|I_{\beta+i}^n|+|I_{\beta+i+1}^n|}
lrd(n,\beta+i) \pm \cO(M_1(n,\beta,p)) \ .
$$
Thus,
\begin{eqnarray*}
lrd(n-1,\alpha)   & \in & lrd(n,\beta)
\frac{|I_{\alpha}^n|+|I_{\alpha+1}^n|}
{|I_{\beta}^n|+|I_{\beta+1}^n|}
\pm   \cO(M_2(n,\beta,p))  \\
&\subset & lrd(n,\beta+i) \frac{|I_{\alpha}^n|+|I_{\alpha+1}^n|}
{|I_{\beta+i}^n|+|I_{\beta+i+1}^n|}
\pm   \cO(M_2(n,\beta,p)) \ .
\end{eqnarray*}
\qed

\subsection{Uniformly asymptotically affine (uaa) homeomorphisms}
\label{sectionXXX1}

The definition of uniformly asymptotically affine  homeomorphism
that we introduce in this paper is more   adapted to our problem
and, apparently, is stronger  than the usual one for symmetric maps,
where   the constant $d$ of the (uua) condition in Definition
\ref{def333}, below,  is taken to be equal to $1$.
However, in Lemma \ref{fewefef}, we will prove that they are equivalent.

\begin{defn}
\label{def333} Let $d\ge 1$ and $\epsilon:\reals^+_0  \to
\reals^+_0$ be a continuous function with $\epsilon (0)=0$. The
homeomorphism $h:I \to  J$ satisfies the {\it $(d, \epsilon)$
uniformly asymptotically affine  condition} if
\begin{equation}
\label{fgggrrs3}
\left|
\log
\frac
{h (x+\delta_2) - h (x)}
{h (x) - h (x-\delta_1)}
\frac
{\delta_1}
{ \delta_2}
\right|
\le \epsilon (\delta_1+\delta_2) \ ,
\end{equation}
for all
$x-\delta_1,x,x+\delta_2 \in I$, such that $\delta_1 > 0$, $\delta_2
>0$ and $ d^{-1} < \delta_2 \ / \delta_1 < d $. The map $h$ is {\it
uniformly asymptotically affine (uaa)} if  for every $d \ge 1$ there
exists   $\epsilon_d$   such that $h$ satisfies the $(d,\epsilon_d)$
uniformly asymptotically affine  condition.
\end{defn}

\begin{lemma}
\label{gfsdgfgd33}
Let   $h:I \to  J$
be a  homeomorphism and $I$ a compact interval.
The following statements are equivalent:
\begin{rlist}
\item
The homeomorphism $h:I \to  J$ is (uaa).
\item
There is a sequence $\gamma_n$ converging to zero,
when $n$ tends to infinity, such that
\begin{equation}
\label{gammaeeffe}
|lrd(n,\beta)| \le \gamma_n \ ,
\end{equation}
for every $n \ge 1$ and every $1 \le \beta < \Omega(n)$.
\end{rlist}
\end{lemma}

\begin{lemma}
\label{fewefef} If  $h:I \to  J$ satisfies the
$(d_0,\epsilon_{d_0})$ uniformly asymptotically affine condition
then the homeomorphism  $h$ is (uaa).
\end{lemma}

\begin{lemma}
\label{lll111}
Let   $h:I \to  J$
be a  homeomorphism and  ${\mathcal G}_\Omega$ a grid of the compact
interval
$I$.
\begin{rlist}
\item
If $h:I \to  J$ is (uaa)
then  there is  a sequence $\alpha_n$ converging to zero,
when $n$ tends to infinity,
		such that
$$
|crd(n,\beta)| \le \alpha_n  \ ,
$$
for every $n \ge 1$ and every $1 \le \beta < \Omega(n)-1$.
\item
If there is a sequence $\alpha_n$ converging to zero, when $n$ tends
to infinity,  such that for every $n \ge 1$ and every $1 \le \beta <
\Omega(n)-1$
$$
|crd(n,\beta)| \le   \alpha_n \ ,
$$
then, for
every  closed interval $K$ contained in the interior of $I$, the
homeomorphism  $h$ is (uaa) in $K$.
\end{rlist}
\end{lemma}

\bigskip
\noindent
{\bf Proof of Lemma \ref{gfsdgfgd33}.}
{\it Let us prove that statement (i) implies statement (ii)}.
Let ${\mathcal G}_\Omega$ be a $(B,M)$ grid of $I$. We have that
\begin{equation}
\label{hhyhhnnbb}
B^{-1} \le r(n,\beta) \le B \ ,
\end{equation}
for every level $n \ge 1$ and every $1 \le \beta < \Omega(n)$.
For every level $n \ge 1$ and every  $1 \le \beta < \Omega(n)$,
let $x-\delta_1, x , x+\delta_2 \in I$     be such that
$I_\beta^n=[x-\delta_1,x]$ and $I_{\beta+1}^n=[x,x+\delta_2]$.
Hence,
$$
\frac{r_h(n,\beta)}{r(n,\beta)}  =
\frac
{h (x+\delta_2) - h (x)} {h (x) - h (x-\delta_1)}
\frac{\delta_1}{\delta_2} \ .
$$
Since  $h:I \to  J$ is
  $(B,\epsilon_{B})$ uniformly asymptotically affine, we get
\begin{equation}
\label{fgfrgnjjnnnbfced}
   lrd(n,\beta)<
\epsilon_{B}(|I_\beta^n|+|I_{\beta+1}^n|)  \ .
\end{equation}
By Remark \ref{gridremark}, there is $B_2=B_2(B,M) < 1$ such that
$|I_\beta^n| \le B_2^n |I|$ and $|I_{\beta+1}^n| \le B_2^n |I|$.
Let $\alpha_n =
\epsilon_{B}(2B_2^n |I|)$.
Hence, by inequality (\ref{fgfrgnjjnnnbfced}), we have
\begin{eqnarray*}
lrd(n,\beta) & < &
\epsilon_{B}(|I_\beta^n|+|I_{\beta+1}^n|) \\
& < & \epsilon_{B}(2B_2^n |I|) \\
& < &  \alpha_n  \ ,
\end{eqnarray*}
for every $n$ and every $1 \le \beta <\Omega(n)$.
Since $\epsilon_{B}(0)=0$ and
$\epsilon_{B}$ is continuous at $0$, we get that
$\alpha_n=\epsilon_{B}(2B_2^n |I|)$ converges to zero, when $n$
tends to infinity.

\bigskip

\noindent {\it Let us prove that statement (ii) implies statement
(i)}. Let ${\mathcal G}_\Omega$ be a $(B,M)$ grid of $I$. Let $d \ge
1$. Let  $x-\delta_1, x, x+\delta_2\in I$, be such that $\delta_1 >
0$, $\delta_2 >0$ and $d^{-1} \le \delta_2 / \delta_1 \le d$. For
every $\alpha > 1$, let $L_1$, $L_2$, $R_1$ and $R_2$ be the
intervals as constructed in Lemma \ref{lll33}. By inequality
(\ref{505050}) and by monotonicity of the homeomorphism $h$, we
obtain that
\begin{equation}
\label{fgrergb33}
\frac{|h(R_1)|}{|h(L_2)|}
\frac{|L_1|}{|R_2|}
\le
\frac
{h (x+\delta_2) - h (x)} {h (x) - h (x-\delta_1)}
\frac
{\delta_1}{ \delta_2}
\le
\frac{|h(R_2)|}{|h(L_1)|}
\frac{|L_2|}{|R_1|} \ .
\end{equation}
By inequality (\ref{tgvcdews}),
\begin{equation}
\label{fgrergb55}
1 \le
\frac{|L_2|}{|L_1|}  \frac{|R_2|} {|R_1|}
\le \alpha^4
\end{equation}
By inequalities (\ref{fgrergb33}) and  (\ref{fgrergb55}), we get
\begin{equation}
\label{fgrergb77}
\alpha^{-4} \frac{|h(R_1)|}{|h(L_2)|}
\frac{|L_2|}{|R_1|}
\le
\frac
{h (x+\delta_2) - h (x)} {h (x) - h (x-\delta_1)}
\frac
{\delta_1}{ \delta_2}
\le  \alpha^4
\frac{|h(R_2)|}{|h(L_1)|}
\frac{|L_1|}{|R_2|}  \ .
\end{equation}
Recalling  equality (\ref{tgvcdews1}) in the proof of Lemma
\ref{lll0}, we have
\begin{eqnarray}
\label{gammaeeffe3353}
\frac{|h(R_2)|}{|h(L_1)|} \frac{|L_1|}{|R_2|} & = &
\frac
{\sum_{i=m-1}^{r-1} \prod_{j=l}^{i} r_h(n_0+n_1,j)  }
{\sum_{i=m-1}^{r-1} \prod_{j=l}^{i} r(n_0+n_1,j)}
\frac
{1+ \sum_{i=l+1}^{m-2} \prod_{j=l}^{i} r(n_0+n_1,j)}
{1+ \sum_{i=l+1}^{m-2} \prod_{j=l}^{i} r_h(n_0+n_1,j) }  \ , \nonumber \\
&&
\\
\frac{|h(R_1)|}{|h(L_2)|} \frac{|L_2|}{|R_1|} & = & \frac
{\sum_{i=m}^{r-2} \prod_{j=l}^{i} r_h(n_0+n_1,j)  }
{\sum_{i=m}^{r-2} \prod_{j=l}^{i} r(n_0+n_1,j)} \frac {1+
\sum_{i=l}^{m-1} \prod_{j=l}^{i} r(n_0+n_1,j)} {1+ \sum_{i=l}^{m-1}
\prod_{j=l}^{i} r_h(n_0+n_1,j) }  \nonumber \ .
\end{eqnarray}
By inequality (\ref{gammaeeffe}), there is $C_0 \ge 1$ and there is a
sequence
$\gamma_n$ converging to zero,
when $n$ tends to infinity,  such that
\begin{equation}
\label{gammaeeffe3}
\frac
{r_h(n_0+n_1,j)}
{r(n_0+n_1,j)}
\in 1 \pm C_0 \gamma_{n_0+n_1}  \ ,
\end{equation}
for every $n_0+n_1$ and for every $1 \le j < \Omega(n_0+n_1)$.
Without loss of generality, we will consider that $\gamma_n$ is a decreasing
sequence.
Hence, by inequalities (\ref{gammaeeffe3353}) and (\ref{gammaeeffe3}),
there is $C_1=C_1(C_0,n_2)  > 1$ such that
$$
\left | \log \frac{|h(R_1)|}{|h(L_2)|} \frac{|L_2|}{|R_1|}   \right |
  \le    C_1 \gamma_{n_0+n_1}   ~~{\rm and}~~
\left | \log \frac{|h(R_2)|}{|h(L_1)|} \frac{|L_1|}{|R_2|}  \right |
  \le
C_1 \gamma_{n_0+n_1} \ .
$$
Therefore, by  inequality (\ref{fgrergb77}), we obtain that
$$
\left| \log   \frac
{h (x+\delta_2) - h (x)} {h (x) - h (x-\delta_1)}
\frac
{\delta_1}{ \delta_2} \right |
\le C_1 \gamma_{n_0+n_1}  + 4 \log(\alpha)  \ .
$$
For every $m=1,2,\ldots$, let $\alpha_m = \exp(1/8m)$.
Hence, we get
\begin{equation}
\label{fgrergb99}
\left| \log   \frac
{h (x+\delta_2) - h (x)} {h (x) - h (x-\delta_1)}
\frac
{\delta_1}{ \delta_2} \right |
\le C_1 \gamma_{n_0+n_1}  + 1/(2m)  \ .
\end{equation}
By Lemma \ref{lll33}, $n_0=n_0( x-\delta_1,x,x+\delta_2) \ge 1$ is
the biggest integer such that $ [x-\delta_1,x+\delta_2] \subset
I_{\beta}^{n_0} \cup I_{\beta+1}^{n_0} $. Hence, there is
$I_\alpha^{n_0+1}  \subset [x-\delta_1,x+\delta_2]$, Thus,
$|I_\alpha^{n_0+1}| \le \delta$, where $\delta=\delta_1+\delta_2$.
By Remark \ref{gridremark}, there is $0 < B_1(B,M) <1$ such that
$|I_\alpha^{n_0+1}| \ge B_1^{n_0+1} |I|$. Hence, we get that

$
B_1^{n_0+1} |I|\le |I_\alpha^{n_0+1}| \le \delta \ ,
$
and so
$$
n_0 \ge \frac{\log \left (\delta  B_1^{-1}
|I|^{-1}\right)}{\log(B_1) }\ .
$$
Therefore,  there is a monotone sequence $\delta_m > 0$ converging
to zero, when $m$ tends to infinity, with the following property: if
$\delta_1+\delta_2 \le \delta_m$ then $n_0=n_0(
x-\delta_1,x,x+\delta_2)$ is sufficiently large such that $C_1
\gamma_{n_0+n_1} \le 1/(2m)$. Hence, by inequality
(\ref{fgrergb99}), for every $m \ge 1$ and every $\delta_0+\delta_1
\le \delta_m$, we have
\begin{equation}
\label{fgrergb93339}
\left| \log   \frac
{h (x+\delta_2) - h (x)} {h (x) - h (x-\delta_1)}
\frac
{\delta_1}{ \delta_2} \right |
\le C_1 \gamma_{n_0+n_1}  + 1/(2m) \le 1/m \ .
\end{equation}
Therefore, we define the continuous function $\epsilon_D:\reals^+ \to
\reals^+$
as follows:
\begin{rlist}
\item $\epsilon_d(\delta_m)=1/(m-1)$ for every $m=2,3,\ldots$;
\item $\epsilon_d$ is affine in every interval $[\delta_m,\delta_m-1]$;
\item Since  $I$ is a compact interval and $h$ is a homoeomorphism,
there is an extension of $\epsilon_d$ to $[\delta_2,\infty)$ such
that inequality (\ref{fgggrrs3}) is satisfied.
\end{rlist}
By inequality (\ref{fgrergb93339}),
we get that   $\epsilon_d$
  satisfies  inequality (\ref{fgggrrs3}).
\qed

\bigskip
\noindent {\bf Proof of Lemma \ref{fewefef}:}  Similarly to the
proof that statement (i) implies statement (ii) of Lemma
\ref{gfsdgfgd33}, we obtain that if $h:I \to  J$ satisfies the
$(d_0,\epsilon_{d_0})$  uniformly asymptotically affine condition
then   satisfies statement (ii) of Lemma \ref{gfsdgfgd33} with
respect to a symmetric grid (see definition of a symmetric grid in
Remark \ref{gridremark}). Since statement (ii) implies statement (i)
of Lemma \ref{lll0}, we get that  the homeomorphism $h$ is (uaa).
\qed

\bigskip

Before proving  Lemma \ref{lll111}, we will
  state and prove Lemma     \ref{reewdsfrewa}
which we will use in the proof of Lemma \ref{lll111}.

\begin{lemma}
\label{reewdsfrewa} Let $h:I \subset \reals \to J \subset \reals$ be
a  homeomorphism and  ${\mathcal G}_\Omega$ a grid of the closed
interval $I$. For every level $n$ and every $0 \le i < \Omega(n)-1$,
let $a(n,i)$ and  $b(n,i)$ be given by
$$ a(n,i)   =
\frac
{1+r_h(n,i)}
{1+r(n,i)}
~~~{\rm and}~~~
b(n,i)  =  \exp(-crd(n,i))  \nonumber \ .
$$
\begin{rlist}
\item
Then, for every $1 \le i <\Omega(n)-1$, we have
\begin{equation}
\label{feress2}
a(n,i) a(n,i-1) b(n,i-1)  =
\frac
{r_h(n,i)}
{r(n,i)}  \ .
\end{equation}
\item
Let  $n \ge 1$ and $\beta,p \in\{2,\ldots,\Omega(n)-1\}$ have the following
properties:
\begin{rlist}
\item [(a)]  There is  $\epsilon >1$ such that $a(n,\beta) \ge \epsilon$.
\item [(b)] There is $\gamma < 1$ such that
$\gamma \le b(n,\beta+i) \le \gamma^{-1}$, for every  $0 \le i < p$.
\end{rlist}
Then,  for every $1 \le i \le p$, we have
\begin{eqnarray}
\label{feress3} a(n,\beta+i) & \ge & 1+
\frac{(\epsilon-1)\gamma^{i}}{2}
\prod_{k=1}^i r(n,\beta+k)  \nonumber \\
&&+\frac{(\epsilon-1)\gamma^{i}B^{-i}}{2} + B (\gamma-1)
\frac{1-(B\gamma^{-1})^{i}}{1-(B\gamma^{-1})}
\ ,
\end{eqnarray}
where $B \ge 1$ is given by the bounded geometry property of the
grid.
\end{rlist}
\end{lemma}

\bigskip
\noindent
{\bf Proof:} {\it  Let us prove inequality (\ref{feress2})}.
By hypotheses, we have
\begin{eqnarray*}
b(n,i-1) & = & \exp(-crd(n,i-1)) \\
& = &
\frac
{1+r(n,i-1)}
{1+r_h(n,i-1)}
\frac
{1+r(n,i)^{-1}}
{1+r_h(n,i)^{-1}} \\
& = &
a(n,i-1)^{-1}
\frac
{1+r(n,i)}
{1+r_h(n,i)}
\frac
{r_h(n,i)}
{r(n,i)} \\
& = &
a(n,i-1)^{-1} a(n,i)^{-1}
\frac
{r_h(n,i)}
{r(n,i)}  \ .
\end{eqnarray*}
Thus,
$$
b(n,i-1) a(n,i-1) a(n,i)    =
\frac
{r_h(n,i)}
{r(n,i)}  \ .
$$

\bigskip
\noindent {\it  Let us prove inequality (\ref{feress3})}. By
definition of  $a(n,i)$ and by  equality  (\ref{feress2}), we have
\begin{eqnarray*}
a(n,i) & = &
\frac
{1+r_h(n,i)}
{1+r(n,i)} \\
b(n,i-1) a(n,i-1) a(n,i)  & = &
\frac
{r_h(n,i)}
{r(n,i)}
\end{eqnarray*}
Hence, we get
\begin{eqnarray*}
a(n,i) (1+r(n,i))& = &
1+r_h(n,i) \\
r_h(n,i)  & = &
b (n,i-1) a(n,i-1) a(n,i)  r(n,i)   \ .
\end{eqnarray*}
Thus,
$$
a(n,i) (1+r(n,i))  =
1+b(n,i-1) a(n,i-1) a(n,i)  r(n,i) \ ,
$$
and so
$$
a(n,i) = (1 - r(n,i)(b(n,i-1)(a(n,i-1)-1)+b(n,i-1)-1
)^{-1} \ .
$$
Therefore, for every  $n \ge 1$, $\beta,p
\in\{2,\ldots,\Omega(n)-1\}$ and $1 \le i \le p$, we get
$$
a(n,\beta+i)  -1 \ge  r(n,\beta+i)(
b(n,\beta+i-1)(a(n,\beta+i-1)-1)+b(n,\beta+i-1)-1)    \ .
$$
Hence, by induction in $1 \le i \le p$, we get
\begin{eqnarray}
\label{fffw0} a(n,\beta+i) -1 & \ge & (a(n,\beta)-1)\prod_{k=1}^i
r(n,\beta+k)
b(n,\beta+k-1)  \nonumber \\
&& + r(n,\beta+i) \sum_{k=1}^i (b(n,\beta+k-1)-1) \prod_{l=k}^{i-1}
r(n,\beta+l)b(n,\beta+l) \ .
\end{eqnarray}
Using that   $B^{-1}< r(n,\beta+k) < B$ by the bounded geometry
property of the grid, we get
\begin{eqnarray}\label{fffw1}
&& (a(n,\beta)-1)   \prod_{k=1}^i   r(n,\beta+k) b(n,\beta+k-1)
   \ge (\epsilon-1)\gamma^{i} \prod_{k=1}^i r(n,\beta+k) \nonumber \\
&   & ~~~~~~~~~~~~~~~~~\ge\frac{(\epsilon-1)\gamma^{i}}{2}
\prod_{k=1}^i r(n,\beta+k) +\frac{(\epsilon-1)\gamma^{i}B^{-i}}{2}
\ .
\end{eqnarray}
Furthermore, noting that $\gamma-1 < 0$, we have
\begin{eqnarray}\label{fffw3}
&& r(n,\beta+i)   \sum_{k=1}^i  (b(n,\beta+k-1)-1) \prod_{l=k}^{i-1}
r(n,\beta+l)b(n,\beta+l)   \\
&&~~~~~~~~~~~~~~~ \ge B (\gamma-1) \sum_{k=1}^i (B
\gamma^{-1})^{i-k} \nonumber \\
&   & ~~~~~~~~~~~~~~~ \ge B (\gamma-1)
\frac{1-(B\gamma^{-1})^{i}}{1-(B\gamma^{-1})} \ .
\end{eqnarray}
Putting inequalities (\ref{fffw0}), (\ref{fffw1}) and (\ref{fffw3})
together, we obtain that
$$
  a(n,\beta+i) -1 \ge \frac{(\epsilon-1)\gamma^{i}}{2}
\prod_{k=1}^i r(n,\beta+k)
+\frac{(\epsilon-1)\gamma^{i}B^{-i}}{2} + B (\gamma-1)
\frac{1-(B\gamma^{-1})^{i}}{1-(B\gamma^{-1})}
\ .
$$
\qed

\bigskip
\noindent
{\bf Proof of Lemma \ref{lll111}.}
{\it Let us prove statement  (i)}.
By   Lemma \ref{gfsdgfgd33}, there is a sequence $\alpha_n$
converging to zero, when $n$ tends to infinity,
		such that
\begin{equation}
\label{fgrgrggggsawqqq3} |lrd(n,\beta)| \le \gamma_n   \ ,
\end{equation}
for every $n \ge 1$ and every $1 \le \beta < \Omega(n)$. By
inequality (\ref{batatas3}) in Lemma \ref{equalities},
we have that
\begin{equation}
\label{fgrgrggggsawqqq1} crd(n,\beta)   \in
\frac{lrd(n,\beta)}{1+r(n,\beta)^{-1}} -
\frac{lrd(n,\beta+1)}{1+r(n,\beta+1)} \pm \cO
(lrd(n,\beta)^2,lrd(n,\beta+1)^2)  \ .
\end{equation}
By the bounded geometry property of a grid, there is $B \ge 1$ such
that $B^{-1} \le r(n,\beta) \le B$. Thus, there is $C_0 > 1$ such
that
\begin{equation}
\label{fgrgrggggsawqqq2} C_0^{-1} \le \frac{1}{1+r(n,\beta)^{-1}}
\le C_0 ~~~{\rm and}~~~ C_0^{-1} \le \frac{1}{1+r(n,\beta+1)} \le
C_0  \ .
\end{equation}
Therefore, putting together
inequalities (\ref{fgrgrggggsawqqq3}),  (\ref{fgrgrggggsawqqq1}) and
(\ref{fgrgrggggsawqqq2}), we obtain that
there is $C_1 > 1$ such that $|crd(n,\beta)| \le C_1 \gamma_n$, for
every level $n$ and every $1 \le \beta < \Omega(n)-1$.

\bigskip

\noindent {\it Let us prove statement (ii)}. Let us  suppose, by
contradiction, that there is $\epsilon_0 > 0$ such that
$|lrd(n(j),\beta(j))| > \epsilon_0$, where $I_{\beta(j)}^{n(j)}
\subset K$ and
  $n(j)$    tends to infinity, when $j$ tends to infinity.
Hence, there is a subsequence $m_j$ such that
  either
$lrd(n(m_j),\beta(m_j)) < - \epsilon_0$ for every $j \ge 1 $, or
$lrd(n(m_j),\beta(m_j))
>\epsilon_0$ for every $j \ge 1 $.
For simplicity of notation, we will denote $n(m_j)$ by $n_j$, and
$\beta(m_j)$ by $\beta_j$.
It is enough to consider the case
where  $lrd(n_j,\beta_j) > \epsilon_0$ (if necessary, after
re-ordering   all the indices). Thus,
there is $\epsilon=\epsilon(\epsilon_0) > 1$ such that, for every
$j \ge 1$,
\begin{equation}
\label{sdfrgeas1}
\frac
{1+r_h(n_j,\beta_j)}
{1+r(n_j,\beta_j)}
>\epsilon \ .
\end{equation}
Let $a(n,i)$ and  $b(n,i)$ be defined as in Lemma \ref{reewdsfrewa}:
\begin{eqnarray}
\label{fghbhhrfs}
a(n,i) & = &
\frac
{1+r_h(n,i)}
{1+r(n,i)}   \\
b(n,i) & = & \exp(-crd(n,\beta))  =  \frac { 1+r(n,i) } { 1+r_h(n,i)
} \frac { 1+r(n,i+1)^{-1} } { 1+r_h(n,i+1)^{-1} } \nonumber \ .
\end{eqnarray}
Hence, we have that $a(n_j,\beta_j) \ge \epsilon$ for every $j \ge
1$. By hypotheses, the cross ratio distortion $crd(n,\beta)$
converges uniformly to zero when $n$ tends to infinity. Thus, there
is an inceasing sequence $\gamma_n$   converging  to one, when $n$
tends to infinity,  such that
\begin{equation}
\label{eeeeefffssss}
\gamma_n \le b(n,i) \le \gamma_n^{-1} \ ,
\end{equation}
for every $1 \le i < \Omega(n)-1$. Let $\eta  =
\min\{(\epsilon-1)/4,1/2\}$. For every $j$ large enough, let $p_j$
be the maximal  integer with the following properties: (i)
$\gamma_{n_j}^{p_j} \ge  \eta \ ;$ (ii) $\gamma_{n_j}^{p_j}
(\epsilon-1)/2 \ge \eta \ ;$ and (iii), letting $B \ge 1$ be as
given by the bounded geometry property of the grid,
$$ \frac{(\epsilon-1)\gamma^{i}B^{-p_j}}{2}  \ge B (1-\gamma)
\frac{1-(B\gamma^{-1})^{p_j}}{1-(B\gamma^{-1})} \ .
$$
Since $\gamma_{n_j}$ converges to one, when $j$ tends to infinity,
we obtain that $p_j$ also tends to infinity, when $j$ tends to
infinity. By properties (ii) and (iii) of $\eta$ and by inequality
(\ref{feress3}), for every $j$ large enough, and for every $1\le i
\le p_j$, we have
\begin{equation}
\label{rfgghggg1} a(n_j,\beta_j+i) \ge 1 + \eta \prod_{k=1}^i
r(n_j,\beta_j+k) > 1 \ .
\end{equation}
For every $j \ge 1$, let $N_j$ be the smallest integer such that
there are four grid intervals $I_{\alpha_j-1}^{N_j}$,
$I_{\alpha_j}^{N_j}$, $I_{\alpha_j+1}^{N_j}$ and
$I_{\alpha_j+2}^{N_j}$ such that
$$
I_{\beta_j}^{n_j} \subset I_{\alpha_j-1}^{N_j}
~~~{\rm and }~~~
I_{\alpha_j}^{N_j}   \cup I_{\alpha_j+1}^{N_j}   \cup I_{\alpha_j+2}^{N_j}
\subset \cup_{i=1}^{p_j-1} I_{\beta_j+i}^{n_j} \ .
$$
Since the grid intervals $I_{\beta_j}^{n_j},
\ldots,I_{\beta_j+p(j)-1}^{n_j}$ are contained in at most four grid
intervals at level $N_j-1$, we obtain that
$$
4M^{n_j-(N_j-1)} \ge p_j \ .
$$
where $M > 1$ is given by the bounded geometry property of the grid.
Thus, $n_j-N_j$ tends to infinity, when $j$ tends to infinity. Let
us denote by $RD(j)$ the following ratio:
\begin{eqnarray*}
\label{rfgghggg3}
RD(j) & = &
\frac
{|I_{\alpha_j}^{N_j}|}
{|J_{\alpha_j}^{N_j}|}
\frac
{|J_{\alpha_j+1}^{N_j}|+|J_{\alpha_j+2}^{N_j}|}
{|I_{\alpha_j+1}^{N_j}|+|I_{\alpha_j+2}^{N_j}|}  \nonumber \\
& = &
\frac
{r_h(N_j,\alpha_j)(1+r_h(N_j,\alpha_j+1))}
{r (N_j,\alpha_j)(1+r (N_j,\alpha_j+1))} \ .
\end{eqnarray*}
By the bounded geometry property of a grid, we have $ B^{-1} <
r(N_j,\alpha_j+i) < B$ for every $-1 \le i \le 3$ and $j \ge 0$. By
Lemma \ref{lll0} and statement (ii) of Lemma \ref{lll1},   there is
$k_0 > 1$ such that $ k_0^{-1} < r_h(N_j,\alpha_j+i) < k_0$ for
every $-1 \le i \le 3$ and $j \ge 0$. Hence, there is $k=k(B,k_0) >
1$ such that for every  $j \ge 0$, we have
\begin{equation}
\label{contradiction}
k^{-1} \le     RD(j) \le   k   \ .
\end{equation}
Now, we are going to prove that
$RD(j)$
tends to infinity,
when $j$ tends to infinity, and so we will get a contradiction.
Let $e_1 < e_2 < e_3< e_4$
be such that
$$
I_{\alpha_j}^{N_j} = \cup_{i=e_1}^{e_2} I_{\beta_j+i}^{n_j}~~~,~~~
I_{\alpha_j+1}^{N_j} = \cup_{i=e_2+1}^{e_3} I_{\beta_j+i}^{n_j}~~~,~~~
I_{\alpha_j+2}^{N_j} = \cup_{i=e_3+1}^{e_4} I_{\beta_j+i}^{n_j} \ .
$$
Hence, we get
\begin{equation}
\label{RDJ} RD(j)     =      \frac {|I_{\alpha_j}^{N_j}|}
{|J_{\alpha_j}^{N_j}|} \frac
{|J_{\alpha_j+1}^{N_j}|+|J_{\alpha_j+2}^{N_j}|}
{|I_{\alpha_j+1}^{N_j}|+|I_{\alpha_j+2}^{N_j}|}
		 =
\frac{R_{1}(j)}{R_{h,1}(j)}
\frac{R_{h,2}(j)+R_{h,3}(j)}{R_{2}(j)+R_{3}(j)} \ ,
\end{equation}
where
\begin{eqnarray*}
R_{1}(j) & = & \frac{|I_{\alpha_j}^{N_j}|} {|I_{\beta_j+e_2}^{n_j}|}
=
1+ \sum_{q=e_1}^{e_2-2}  \prod_{i=q+1}^{e_2-1} r(n_j,\beta_j+i)^{-1}  \\
R_{h,1}(j) & = & \frac {|J_{\alpha_j}^{N_j}|}
{|J_{\beta_j+e_2}^{n_j}|}  =
1+ \sum_{q=e_1}^{e_2-2}  \prod_{i=q+1}^{e_2-1} r_h(n_j,\beta_j+i)^{-1}  \\
R_{2}(j) &= &
\frac{|I_{\alpha_j+1}^{N_j}|}
{|I_{\beta_j+e_2}^{n_j}|}
=
\sum_{q=e_2}^{e_3-1}  \prod_{i=e_2}^{q} r  (n_j,\beta_j+i) \\
R_{h,2}(j) &  = &
\frac{|J_{\alpha_j+1}^{N_j}|}
{|J_{\beta_j+e_2}^{n_j}|}
=  \sum_{q=e_2}^{e_3-1}  \prod_{i=e_2}^{q} r_h  (n_j,\beta_j+i) \\
R_{3}(j) &= &
\frac{|I_{\alpha_j+2}^{N_j}|}
{|I_{\beta_j+e_2}^{n_j}|}
=
\sum_{q=e_3}^{e_4-1}  \prod_{i=e_2}^{q} r  (n_j,\beta_j+i) \\
R_{h,3}(j) &  = &
\frac{|J_{\alpha_j+2}^{N_j}|}
{|J_{\beta_j+e_2}^{n_j}|}
=
\sum_{q=e_3}^{e_4-1}  \prod_{i=e_2}^{q} r_h (n_j,\beta_j+i) \ .
\end{eqnarray*}
Hence, by inequalities (\ref{fghbhhrfs}) and (\ref{rfgghggg1}), for
every $1\le i \le p_j$, we get
\begin{equation}
\label{rfgghggg2} \frac {r_h(n_j,\beta_j+i)} {r(n_j,\beta_j+i)} > 1
\ .
\end{equation}
Thus, we deduce that
\begin{eqnarray}
\label{R1} R_{h,1}(j) & = & 1+ \sum_{q=e_1}^{e_2-2}
\prod_{i=q+1}^{e_2-1} r(n_j,\beta_j+i)^{-1} \frac {r(n_j,\beta_j+i)}
{r_h(n_j,\beta_j+i)}  \nonumber\\
&\le &  1+ \sum_{q=e_1}^{e_2-2}  \prod_{i=q+1}^{e_2-1} r(n_j,\beta_j+i)^{-1}
\nonumber\\
& = &  R_1(j) \ .
\end{eqnarray}
By inequality (\ref{rfgghggg2}), we obtain
\begin{eqnarray}
\label{R2}
R_{h,2}(j)
& =  &  \sum_{q=e_2}^{e_3-1}  \prod_{i=e_2}^{q} r  (n_j,\beta_j+i)
\frac{r_h  (n_j,\beta_j+i)}
{r   (n_j,\beta_j+i)} \nonumber \\
  & \ge  &  R_{2}(j) \ .
\end{eqnarray}
Now, let us bound $R_{h,3}(j)$ in terms of $R_{3}(j)$.
Putting together  inequalities (\ref{feress2}) and (\ref{rfgghggg1}), we
obtain
\begin{eqnarray}
\label{grdreeee}
\frac
{r_h(n,i)}
{r(n,i)}
& = &
b(n,i-1) a(n,i) a(n,i-1) \nonumber\\
& \ge & b(n,i-1) a(n,i) \  .
\end{eqnarray}
Noting that $e_3-e_2 < p_j$, and by
inequality (\ref{eeeeefffssss}) and  property (i) of $\eta$, we get
$$
\prod_{i=e_2}^{e_3-1} b(n_j,\beta_j+i-1) \ge  \gamma^{p_j} \ge \eta
\ .
$$
Hence, by inequalities (\ref{rfgghggg1}) and (\ref{grdreeee}), we get
\begin{eqnarray}
\label{rfgergfcdff1} \prod_{i=e_2}^{e_3-1} \frac
{r_h(n_j,\beta_j+i)} {r(n_j,\beta_j+i)} & \ge &
\prod_{i=e_2}^{e_3-1} b(n_j,\beta_j+i-1)
a(n_j,\beta_j+i)  \nonumber  \\
  & \ge &
\eta \prod_{i=e_2}^{e_3-1}
\left(1+   \eta \prod_{k=1}^i r(n_{m_i},\beta_j+k) \right) \nonumber  \\
& \ge &
\eta \left(1+  \eta   \sum_{i=e_2}^{e_3-1} \prod_{k=1}^i
r(n_{m_i},\beta_j+k)
\right) \nonumber \\
& \ge &
   \eta^2 \frac{   |I_{\alpha_j+1}^{N_j}|}
{|I_{\beta_j+1}^{n_j}|}    \ .
\end{eqnarray}
Noting that $I_{\beta_j+1}^{n_j} \subset I_{\alpha_j-1}^{N_j} \cup
I_{\alpha_j}^{N_j}$ and by  the  bounded geometry property of the
grid, we get
  \begin{equation}
\label{rfgergfcdff2}
\frac{|I_{\alpha_j+1}^{N_j}|}{|I_{\beta_j+1}^{n_j}|}
  \ge   B^{-2} B_2^{N_j-n_j} \ ,
\end{equation}
where $B_2 < 1$ is given in Remark \ref{gridremark}.
Putting together inequalities (\ref{rfgergfcdff1}) and
(\ref{rfgergfcdff2}),  we obtain that
$$
\prod_{i=e_2}^{e_3-1} \frac {r_h(n_j,\beta_j+i)} {r(n_j,\beta_j+i)}
  \ge
\eta^2 B^{-2} B_2^{N_j-n_j} \ .
$$
Hence,
\begin{eqnarray}
\label{R3} R_{h,3}(j) &= &  \prod_{i=e_2}^{e_3-1} \frac{r_h
(n_j,\beta_j+i)}{r (n_j,\beta_j+i)} r (n_j,\beta_j+i)
\sum_{q=e_3}^{e_4-1}  \prod_{i=e_3}^{q} \frac{r_h (n_j,\beta_j+i)}{r
(n_j,\beta_j+i)}
r(n_j,\beta_j+i) \nonumber \\
&\ge & \eta^2 B^{-2} B_2^{N_j-n_j} \prod_{i=e_2}^{e_3-1} r
(n_j,\beta_j+i) \sum_{q=e_3}^{e_4-1}  \prod_{i=e_3}^{q}
r(n_j,\beta_j+i) \nonumber \\
&= & \eta^2 B^{-2} B_2^{N_j-n_j} R_{3}(j) \ .
\end{eqnarray}
Noting that $R_{2}(j)  R_{3}(j)^{-1} =|I_{\alpha_j+1}^{N_j}|
|I_{\alpha_j+2}^{N_j}|^{-1}$ and by the bounded geometry property of
the grid, we obtain
$$B^{-1} \le R_{2}(j)  R_{3}(j) ^{-1} \le B  \ .$$  Therefore, putting
together inequalities (\ref{R1}), (\ref{R2}) and (\ref{R3}), we
obtain that
\begin{eqnarray*}
RD(j)
& = &
\frac{R_{1}(j)}{R_{h,1}(j)}
\frac{R_{h,2}(j)+R_{h,3}(j)}{R_{2}(j)+R_{3}(j)}  \\
& \ge & \frac{R_{2}(j)+  \eta^2 B^{-2} B_2^{N_j-n_j} R_{3}(j)
}{R_{2}(j)+R_{3}(j)} \\
& \ge & \frac{1+\eta^2 B^{-3} B_2^{N_j-n_j} }{1+B} \ .
\end{eqnarray*}
Since $B_2^{N_j-n_j}$ tends to infinity, when $j$ tends to infinity,
we get that $RD(j)$ also tends to infinity, when $j$ tends to
infinity. However, by inequality (\ref{contradiction}), this is
absurd. \qed

\subsection{$C^{1+r}$ diffeomorphisms}
\label{hhholderrr}

Let $0 < r \le 1$. We say that  a  homeomorphism   $h:I \to  J$
is $C^{1+r}$ if its differentiable and its first derivative $dh:I \to
\reals$
is $r$-H\"older continuous, i.e. there is $C \ge 0$ such that,
for every $x,y \in I$,
$$
|dh(y)-dh(x)| \le C |y-x|^r \ .
$$
In particular, if $r=1$ then $dh$ is  Lipschitz.

\begin{lemma}
\label{r1+}
   Let   $h:I \to  J$
be a  homeomorphism, and let $I$ be  a compact interval with a grid
${\mathcal G}_\Omega$.
\begin{rlist}
\item  For $0 < r \le 1$, the map $h$ is a $C^{1+r}$ diffeomorphism if,
and only if,
for every $n \ge 1$ and for every $1 \le \beta < \Omega(n)$,
we have that
\begin{equation}
\label{g1+}
|lrd(n,\beta)| \le \cO(|I_\beta^n|^{r}) \ .
\end{equation}
\item
The map $h$ is affine if, and only if,
for every $n \ge 1$ and every $1 \le \beta < \Omega(n)$,
we have that
\begin{equation}
\label{g1+++}
|lrd(n,\beta)| \le \co(|I_\beta^n|) \ .
\end{equation}
\end{rlist}
\end{lemma}

\begin{lemma}
\label{lll11179}
Let $0<r \le 1$. Let   $h:I \to  J$
be a  homeomorphism and
${\mathcal G}_\Omega$ a grid of the compact interval $I$.
\begin{rlist}
\item
If $h:I \to  J$ is  a $C^{1+r}$ diffeomorphism
then,  for every $n \ge 1$ and every $1 \le \beta < \Omega(n)-1$,
we have that
\begin{equation}
\label{78new}
|crd(n,\beta)| \le \cO(|I_\beta^n|^{r}) \ .
\end{equation}
\item
If, for every $n \ge 1$ and every $1 \le \beta < \Omega(n)-1$,
we have that
\begin{equation}
\label{79new}
|crd(n,\beta)| \le \cO(|I_\beta^n|^{r}) \ ,
\end{equation}
then,  for every  closed interval $K$ contained in the interior of $I$,
the homeomorphism $h|K$ restricted to $K$ is  a $C^{1+r}$
diffeomorphism.
\end{rlist}
\end{lemma}

\bigskip
\noindent
{\bf Proof of Lemma \ref{r1+}:}
By the Mean Value Theorem,
if $h$ is a $C^{1+r}$ diffeomorphism  then
for every $n \ge 1$ and for every
grid interval $I_\beta^n$ we get that
$lrd(n,\beta) \in \pm \cO(|I_\beta^n|^r)$, and so inequality (\ref{g1+}) is
satisfied.
If $h$ is affine then,
for every $n \ge 1$ and for every
grid interval $I_\beta^n$, we get that
$lrd(n,\beta) =0$, and so inequality (\ref{g1+++}) is satisfied.

\bigskip
\noindent {\it Let us prove that inequality (\ref{g1+}) implies that
$h$ is $C^{1+r}$}. For every point $P \in I$, let
$I_{\alpha_1}^1,I_{\alpha_2}^2,\ldots$ be a sequence of grid
intervals $I_{\alpha_n}^n$ such that $P \in I_{\alpha_{n}}^{n}$ and
$I_{\alpha_n}^n \subset I_{\alpha_{n-1}}^{n-1}$ for every $n > 1$.
Let us suppose that $I_{\alpha_{n-1}}^{n-1}=\cup_{i=0}{j}
I_{\alpha_n+i}^n$ for some $j=j(\alpha_n) \ge 1$. By inequality
(\ref{g1+})and using the bounded geometry of the grid, we obtain
that
\begin{eqnarray*}
\label{1eeaaae}
   \frac {dh(n-1,\alpha_{n-1})} {dh(n,\alpha_n)}
& = &
\frac {1+ \sum_{i=1}^{j}    \prod_{k=1}^{i} r_h (n,\alpha_n+i) }
{1+ \sum_{i=1}^{j}    \prod_{k=1}^{i} r(n,\alpha_n+i)} \\
& \in &
\frac {1+ \sum_{i=1}^{j}    \prod_{k=1}^{i} r  (n,\alpha_n+i)(1 \pm
\cO(|I_{\alpha_n+i}^n|)}
{1+ \sum_{i=1}^{j}    \prod_{k=1}^{i} r(n\alpha_n+i)} \\
& \subset &    \cO(|I_{\alpha_n}^n|^{r}) \ .
\end{eqnarray*}
A similar argument to the one above implies that for all
$I_{\alpha_n}^n \subset I_{\alpha_{n-1}}^{n-1}$, we have
$$
dh(n,{\alpha_{n}})
  \in
dh(n-1,\alpha_{n-1}) \pm \cO(|I_{\alpha_{n-1}}^{n-1}|^{r}) \ .
$$
Hence, using the  bounded geometry property of a
grid,
for every $m \ge 1$ and for every $n \ge m$, we get
\begin{equation}
\label{AC7}
dh(n,{\alpha_{n}})
  \in
dh(m,\alpha_m)
\pm \cO(|I_{\alpha_m}^m|^{r}) \ .
\end{equation}
Thus, the average derivative $dh(n,{\alpha_{n}})$ converges to a
value $d_P$, when $n$ tends to infinity. Let us prove that $h$ is
differentiable at $P$ and that $dh(P)=d_P$. Let $L$ be any interval
such that the point $P \in L$. Take the largest $m \ge 1$ such that
there is a grid interval $I_\gamma^m$ with the property that $L
\subset   \cup_{j=-1,0,1} I_{\gamma+j}^m$. By the bounded geometry
property of a grid, there is $C \ge 1$, not depending upon $P$, $L$
and $I_\gamma^m$, such that
\begin{equation}
\label{AC69}
C^{-1} \le
\frac{|I_\gamma^m|}{|L|}
\le C \ .
\end{equation}
Then, using inequality (\ref{g1+}) and   the bounded geometry of the
grid, for every $j=\{-1,0,1\}$, we obtain that $$
\left|ldh(m,\gamma+j)-ldh(m,\gamma)\right| \le \cO(|L|^r) \ ,
$$
and so
\begin{equation}
\label{AC6966}   dh(m,\gamma+j) \in dh(m,\gamma) \pm \cO(|L|^r) \ .
\end{equation}
For every $n \ge m$, take the smallest sequence of adjacent grid
intervals $I_{\beta_n}^n, \ldots, I_{\beta_n+i_n}^n$,
at level $n$,
such that $L \subset \cup_{i=0}^{i_n} I_{\beta_n+i}^n \subset
\cup_{j=-1,0,1} I_{\gamma+j}^m$. By inequalities (\ref{AC7}) and
(\ref{AC6966}), for every $I_{\beta_n+i}^n \subset
I_{\gamma+j(i)}^m$ we get that
\begin{eqnarray*}
dh(m,\beta_n+i) & \in &  dh(m,\gamma+j(i)) \pm \cO(|I_{\gamma+j(i)}^m|^r)\\
& \subset & d_P \pm \cO(|L|^r) \ . \end{eqnarray*} Hence,
\begin{eqnarray}
\label{AC9} \label {AC71}
\frac{h(L)}{L} & = &
\lim_{n \to \infty} \sum_{i=0}^{i_n}
\frac{|I_{\beta_n+i}^n|}{|L|}
dh(n,\beta_n+i) \nonumber \\
& \in &
\lim_{n \to \infty} \sum_{i=0}^{i_n}
\frac{|I_{\beta_n+i}^n|}{|L|}
\left( d_P \pm
\cO(|L|^{r}) \right)  \nonumber  \\
& \subset &
d_P \pm
\cO(|L|^{r})
\ .
\end{eqnarray}
Therefore, for every $P \in I$, the homeomorphism $h$ is
differentiable at $P$ and $dh(P)=d_P$. Let us check that $dh$ is
$r$-H\"older continuous. For every $P,P' \in I$, let $L$ be the
closed interval $[P,P']$. Using inequality (\ref{AC71}), we obtain
that
\begin{eqnarray*}
\label{AC73}
dh(P')-dh(P) &\in &
\frac{h(L)}{L} - \frac{h(L)}{L}
\pm
\cO(|L|^{r}) \nonumber \\
&\subset & \pm \cO(|L|^{r}) \ ,
\end{eqnarray*}
and so $dh$ is $r$-H\"older
continuous.

\bigskip
\noindent
{\it Let us prove that inequality (\ref{g1+++}) implies that $h$ is affine.}
A similar argument to the one above
gives us that $h$ is differentiable and
that
\begin{equation}
\label{AC75}
|dh(P')-dh(P)| \le
\co(|P'-P|) \ ,
\end{equation}
for every $P,P' \in I$.
Hence,   we get that
\begin{eqnarray*}
|dh(P')-dh(P)|  & \le & \lim_{n \to \infty} \sum_{i=0}^{n-1}
\left|dh\left(P + \frac{(i+1)(P'-P)}{n} \right) -dh\left(P +
\frac{i (P'-P)}{n}\right)\right| \nonumber \\ & \le & \lim_{n \to
\infty}
n \  \co\left(\frac{P'-P}{n} \right) = 0 \ ,
\end{eqnarray*}
and so $h$ is an affine map.
\qed

\bigskip
\noindent
{\bf Proof of Lemma \ref{lll11179}:}
{\it Proof of statement (i)}:
By Lemma \ref{r1+}, for every $n \ge 1$ and for every $1 \le \beta <
\Omega(n)$,
we have that $|lrd(n,\beta)| \le \cO(|I_\beta^n|^{r})$.
Hence, by the  bounded geometry property of a grid and by   inequality
(\ref{batatas3}),
we get $|crd(n,\beta)| \le \cO(|I_\beta^n|^{r})$.

\bigskip
\noindent
{\it Proof of statement (ii)}:
Let $K$ be a closed interval contained in the interior of $I$.
By Lemmas \ref{gfsdgfgd33} and \ref{lll111}, there is a decreasing sequence
of
positive reals $\epsilon_n$ which
converges to
$0$, when $n$ tends to $\infty$,
such that
\begin{equation}
\label{ffffrrrs} |lrd(n,\beta)| < |\epsilon_n| \ ,
\end{equation}
for all $n \ge 1$ and for all grid interval
  $I_\beta^n$ intersecting $K$.
For every    grid interval $I_{\alpha}^{n-1}$  intersecting
  $K$, let $k_1=k_1(n,\alpha)$ and $k_2=k_2(n,\alpha)$
be such that
$
\cup_{\beta=k_1}^{k_2} I_{\beta}^n = I_{\alpha}^{n-1} \cup
I_{\alpha+1}^{n-1}$. Let the integers $\beta$ and $i$ be such that
$k_1\le \beta \le k_2$ and   $k_1\le \beta +i \le k_2$. By the
bounded geometry property of a grid, and by inequalities
(\ref{batatas4}) and (\ref{78new}), we get
$$ lrd(n,\beta+i)  \in
\pm \cO \left( |lrd(n,\beta)|+(|I_{\beta}^n|+|I_{\beta+1}^n|)^r
\right) \ . $$
Therefore,
\begin{equation}
\label{batatas33}
L_2(n,\beta,p)  \in
\pm \cO \left( lrd(n,\beta)^2+(|I_{\beta}^n|+|I_{\beta+1}^n|)^{2r}
\right)
\ .
\end{equation}
By   inequalities (\ref{batatas5}) and (\ref{batatas33}),
we get
\begin{equation}
\label{bat37}
lrd(n-1,\alpha)   \in
\frac{|I_{\alpha}^{n-1}|+|I_{\alpha+1}^{n-1}|}{|I_{\beta+i}^n|+|I_{\beta+i+1}^n|}
lrd(n,\beta+i)  \pm \cO \left(  lrd(n,\beta)^2 +
(|I_{\beta}^n|+|I_{\beta+1}^n|)^r  \right)  \ .
\end{equation}
Let us suppose, by contradiction, that there is   a sequence of grid
intervals $I_{\beta_j}^{n_j}$ and a sequence of  positive reals
$|e_j|$ which tends to infinity, when $j$ tends to infinity,
such that
\begin{equation}
\label{D1}
lrd(n_j,\beta_j)  =  e_j(|I_{\beta_j}^{n_j}|+|I_{\beta_j+1}^{n_j}|) ^r \ .
\end{equation}
Using that the number of grid intervals at every level $n$ is finite,
we obtain that
there exists a subsequence $m_j$ of $j$ such that
$I_{\beta_{m_{j+1}}}^{n_{m_{j+1}}} \subset I_{\beta_{m_j}}^{n_{m_j}}$.
Therefore, there exists a sequence of grid intervals
$I_{\alpha_1}^1,I_{\alpha_2}^2,\ldots$   with the following
properties:
\begin{rlist}
\item for every $i \ge 1$,
  $I_{\alpha_{i+1}}^{i+1} \subset I_{\alpha_i}^i$;
\item for every $i \ge 1$, let  $a_i$ be determined such that
\begin{equation}
\label{D3}
lrd(i,\alpha_i)   =  a_i
(|I_{\alpha_i}^i|+|I_{\alpha_i+1}^i|)^r \ .
\end{equation}
Then, there is a subsequence $m_j$ of $j$ such that
  $|a_i| \le |a_{m_j}|$ for every  $1 \le i \le m_j$, and $|a_{m_j}|$
  tends to infinity, when $j$ tends to infinity.
\end{rlist}
Let us denote $|I_{\alpha_i}^i|+|I_{\alpha_i+1}^i|$ by $B_i$.
Using inequality (\ref{bat37}) inductively, we get
\begin{eqnarray}
\label{bat3733} lrd(m_j,\beta_{m_j})  & \in & \frac{B_{m_j}}{B_1}
lrd(1,\alpha_1)
  \pm
\cO \left ( \sum_{i=2}^{m_j} \frac{B_{m_j}}{B_i} (lrd(i,\alpha_i)^2
+  B_i^{r}) \right) \ .
\end{eqnarray}
By the  bounded geometry property of a grid, there is  $0<\theta < 1$ such
that
\begin{equation}
\label{neeen}
\frac{B_k}{B_i} \le \theta^{k-i} \ ,
\end{equation}
for every $1 \le i \le m_j$ and for every $1 \le k \le m_j$. Noting
that $|a_1| \le |a_{m_j}|$, by inequalities (\ref{D3}) and
(\ref{neeen}), we get
\begin{eqnarray}
\label{ffrf1}
\frac{B_{m_j}}{B_1} lrd(1,\alpha_1)  & = &   \frac{a_1 B_1^r B_{m_j}}{B_1}
\nonumber \\
& \in & \pm \cO \left (|a_{m_j}| B_{m_j}^r \theta^{(1-r)m_j} \right)
\end{eqnarray}
By inequality (\ref{ffffrrrs}), $a_i B_i \le \epsilon_i$, and $|a_i|
\le |a_{m_j}|$ for $i \le m_j$. Hence, by inequalities (\ref{D3})
and (\ref{neeen}), we obtain that
\begin{eqnarray}
\label{ffrf3} \frac{B_{m_j}}{B_i} (lrd(i,\alpha_i)^2+B_i^{r}) & = &
\frac{a_i (a_i B_i^{r})(B_i^{r}B_{m_j}) +B_i^r
B_{m_j}}{B_i}\nonumber
\\ & \in & \pm \cO \left( (|a_{m_j}| \epsilon_i +1)
B_{m_j}^r  \theta^{(1-r)(m_j-i)}   \right)
\end{eqnarray}
Using inequalities (\ref{ffrf1}) and (\ref{ffrf3}) in inequality
(\ref{bat3733}),
we get
\begin{eqnarray}
\label{bat37335} \frac{|lrd(m_j,\beta_{m_j})|}{|a_{m_j}| B_{m_j}^r}
&\le &    \cO \left (     \theta^{(1-r)m_j} + \sum_{i=2}^{m_j}
\left( (  \epsilon_i + |a_{m_j}|^{-1}) \theta^{(1-r)(m_j-i)}
  \right)
\right)  \nonumber \\
& \le &    \cO \left (     \theta^{(1-r)m_j} +
\frac{|a_{m_j}|^{-1}}{1-\theta^{1-r}}+
\sum_{i=2}^{m_j} \left(   \epsilon_i \theta^{(1-r)(m_j-i)}
  \right)
\right)  \ .
\end{eqnarray}
Since $\epsilon_i$ converges to zero, when $i$ tends to infinity,
inequality  (\ref{bat37335}) implies that  there is $j_0 \ge 0$ such
that, for every $j \ge j_0$, we get
$$
|lrd(m_j,\beta_{m_j})  |  < |a_{m_j}| B_{m_j}^r \ ,
$$
  which contradicts (\ref{D3}).
\qed

\subsection{$C^{2+r}$ diffeomorphisms}

Let $0 < r \le 1$. We say that  a  homeomorphism   $h:I \to  J$
is $C^{2+r}$ if its twice differentiable and its second derivative $d^2h:I
\to
\reals$
is $r$-H\"older continuous.

\begin{lemma}
\label{lll111r2}
Let $0<r \le 1$. Let   $h:I \to  J$
be a  homeomorphism and
${\mathcal G}_\Omega$ a grid of the compact interval $I$.
\begin{rlist}
\item
If $h:I \to  J$ is  $C^{2+r}$
then
$$
|crd(n,\beta)| \le \cO(|I_\beta^n|^{1+r}) \ ,
$$
for every $n \ge 1$ and every $1 \le \beta < \Omega(n)-1$.
\item
If, for every $n \ge 1$ and every $1 \le \beta < \Omega(n)-1$,
we have that
\begin{equation}
\label{C333eee}
|crd(n,\beta)| \le \cO(|I_\beta^n|^{1+r}) \ ,
\end{equation}
then,   for every  closed interval $K$ contained in the interior of $I$,
the homeomorphism $h|K$ restricted to $K$ is  $C^{2+r}$.
\end{rlist}
\end{lemma}

Before proving  Lemma  \ref{lll111r2},
we will
  state and prove Lemma     \ref{equalities333}
which we will use later in the proof of Lemma \ref{lll111r2}.

\begin{lemma}
\label{equalities333} Let ${\mathcal G}_\Omega$ be a grid of the
closed interval $I$. Let $h:I \subset \reals \to J \subset \reals$
be a homeomorphism such that for every $n \ge 1$ and every $1 \le
\beta < \Omega(n)-1$,
\begin{equation}
\label{sdfgeeegggg} |crd(n,\beta)| \le \cO(|I_\beta^n|^{1+r}) \ ,
\end{equation}
where $0 \le r < 1$. Then, for every closed interval $K$ contained
in the interior of $I$, the logarithmic ratio distortion and the
cross ratio distortion satisfy the following estimates:
\begin{rlist}
\item There is a constant $C(i)>0$,   not depending upon the level $n$ and
  not depending  upon $1 \le \beta \le \Omega(n)$, such that
\begin{eqnarray}
\label{batatas43} lrd(n,\beta+i) & \in &
\frac{|I_{\beta+i}^n|+|I_{\beta+i+1}^n|}{|I_{\beta}^n|+|I_{\beta+1}^n|}
lrd(n,\beta) \pm   C(i) |I_{\beta}^n|^{1+r} \ .
\end{eqnarray}
\item
Let $I_{\alpha}^{n-1}$ and $I_{\alpha+1}^{n-1}$ be two adjacent grid
intervals. Let $I_{\beta}^n$ and $I_{\beta+1}^n$ be   grid intervals
contained in the union $I_{\alpha}^{n-1} \cup I_{\alpha+1}^{n-1}$.
Then,
\begin{equation}
\label{batatas53} lrd(n-1,\alpha)   \in
\frac{|I_{\alpha}^{n-1}|+|I_{\alpha+1}^{n-1}|}{|I_{\beta}^n|+|I_{\beta+1}^n|}
lrd(n,\beta) \pm \cO(|I_{\beta}^n|^{1+r}) \ .
\end{equation}
\end{rlist}
\end{lemma}

\bigskip
\noindent {\bf Proof of Lemma \ref{equalities333}:} By Lemma
\ref{lll11179}, for every $0<s < 1$, the homeomorphism $h|K$ is
$C^{1+s}$, and so the map $\psi:I \to \reals$ is well-defined  by
$\psi(x)=\log dh(x)$. By  bounded geometry property of a grid and by
inequality (\ref{sdfgeeegggg}), for every integer $i$, there is a
positive constant $E_1(i)$ such that
\begin{equation}
\label{BAC15} |crd(n,\beta+j_1)| \le E_1(i) (|I_{\beta}^n|^{1+r}) \
,
\end{equation}
for every grid interval $I_{\beta}^n$ and $0 \le j_1 \le i$. Take
$s<1$ such that $2s=1+r$ and  $0 \le j_2 \le i$. By inequality
(\ref{sdfgeeegggg}) and statement (ii) of Lemma \ref{lll11179}, $h$
is $C^{1+s}$. Hence, using   the  bounded geometry property of a
grid and  statement (i) of Lemma \ref{r1+}, we obtain that
\begin{eqnarray}
\label{BACCC13} |lrd(n,\beta+j_1) lrd(n,\beta+j_2)| & \le &
\cO(|I_{\beta+j_1}^n|^{s} |I_{\beta+j_2}^n|^{s}) \nonumber \\
& \le & E_2(i) (|I_{\beta}^n|^{1+r})
\end{eqnarray}
where $E_2(i)$ is a positive constant depending upon $i$. Using
inequalities (\ref{BAC15}) and (\ref{BACCC13}) in (\ref{batatas4}),
we get inequality (\ref{batatas43}). Furthermore, using inequalities
(\ref{BAC15}) and (\ref{BACCC13}) in (\ref{batatas5}), we get
inequality (\ref{batatas53}). \qed

\bigskip
\noindent {\bf Proof of Lemma \ref{lll111r2}:} {\it Proof of
statement (i)}: Let $h$ be $C^{2+r}$ and let  $\psi:I \to \reals$
be given by $\psi(x)=\log dh(x)$.
For every $n \ge 1$, let
$I_\gamma^n=[x,y]$, $I_{\gamma+1}^n=[y,z]$ and $I_{\gamma+1}^n=[z,w]$ be
adjacent grid intervals,
  at     level $n$.
By Taylor series, we get
\begin{eqnarray*}
|h(I_\gamma ^n)| & \in & |I_\gamma ^n| dh(y) + |I_\gamma ^n|^2
d^2h(y) \pm \cO(|I_\gamma ^n|^{2+r}) \\
| h(I_{\gamma +1} ^n)| & \in &
|I_{\gamma +1} ^n| dh(y) - |I_{\gamma +1} ^n|^2 d^2h(y) \pm
\cO(|I_{\gamma +1} ^n|^{2+r}) \\
| h(I_{\gamma +1} ^n)| & \in &
|I_{\gamma +1} ^n| dh(z) + |I_{\gamma +1} ^n|^2 d^2h(z) \pm
\cO(|I_{\gamma +1} ^n|^{2+r}) \\
| h(I_{\gamma +2} ^n)| & \in &
|I_{\gamma +2} ^n| dh(z) - |I_{\gamma +2} ^n|^2 d^2h(z) \pm
\cO(|I_{\gamma +2} ^n|^{2+r}) \ .
\end{eqnarray*}
Therefore,
\begin{eqnarray*}
\frac{|h(I_{\gamma +1} ^n)|}{|I_{\gamma +1} ^n|} \frac{|I_{\gamma}
^n|}{|h(I_{\gamma } ^n)|}
& \in &   \frac{dh(y) - |I_{\gamma +1} ^n| d^2h(y)
\pm \cO(|I_{\gamma +1} ^n|^{1+r})} {dh(y) + |I_{\gamma} ^n|
d^2h(y) \pm \cO(|I_{\gamma} ^n|^{1+r})} \\ & \in &
1- (|I_\gamma^n|+ |I_{\gamma+1}^n|)
\frac {d \psi(y)} {2} \pm \cO((|I_\gamma^n|+ |I_{\gamma+1}^n|)^r) \
, \end{eqnarray*} and so $$lrd(n,\gamma) \in  -(|I_\gamma^n|+
|I_{\gamma+1}^n|) \frac {d \psi(y)} {2} \pm \cO((|I_\gamma^n|+
|I_{\gamma+1}^n|)^r) \ .$$
Similarly, we get
  $$ lrd(n,\gamma+1) \in -(|I_{\gamma+1}^n|+
|I_{\gamma+2}^n|) \frac {d\psi(z)} {2} \pm \cO((|I_{\gamma+1}^n|+
|I_{\gamma+2}^n|)^r). $$ Therefore, by inequality
(\ref{batatas3}), the cross ratio distortion $c(n,\gamma) \in \pm
\cO(|I_\gamma^n|^r)$.

\bigskip
\noindent {\it Proof  of statement  (ii)}: We prove statement  (ii),
first  in the case where $0 < r <1$ and secondly in the case where $r=1$.

\bigskip
\noindent
{\it  Case $0 < r <1$}: By
Lemma \ref{lll11179}, for every $0<s < 1$, the homeomorphism $h|K$
is $C^{1+s}$, and so the map $\psi:I \to \reals$ is well-defined
by $\psi(x)=\log dh(x)$. For every point $P \in I$, let
$I_{\alpha_1}^1,I_{\alpha_2}^2,\ldots$ be a sequence of grid
intervals $I_{\alpha_n}^n$ such that $P \in I_{\alpha_{n}}^{n}$
and $I_{\alpha_n}^n \subset I_{\alpha_{n-1}}^{n-1}$ for every $n >
1$. By the  bounded geometry property of a grid and by inequality
(\ref{C333eee}), for every grid interval $I_{\beta}^n \subset
\cup_{i=-1,0,1} I_{\alpha_{n-1}+i}^{n-1}$, we have that
\begin{equation}
\label{BAC1}
|crd(n,\beta)| \le \cO(|I_{\alpha_n}^n|^{1+r}) \ .
\end{equation}
By inequality  (\ref{batatas53}),
we have
$$
\frac{lrd(n-1,\alpha_{n-1})}{|I_{\alpha_{n-1}}^{n-1}|+|I_{\alpha_{n-1}+1}^{n-1}|}
\in
\frac{lrd(n,\alpha_n)}{|I_{\alpha_n}^n|+|I_{\alpha_n+1}^n|}
   \pm \cO(|I_{\alpha_{n}}^{n}|^r) \ .
$$
Hence, by the  bounded geometry property of a grid,
for every $m \ge 1$ and for every $n \ge m$, we get that
\begin{equation}
\label{BAC7}
\frac{lrd(n,\alpha_{n})}{|I_{\alpha_n}^{n}|+|I_{\alpha_n+1}^{n}|}
  \in
\frac{lrd(m,\alpha_m)}{|I_{\alpha_m}^{m}|+|I_{\alpha_m+1}^{m}|}
\pm \cO(|I_{\alpha_m}^m|^{r}) \ .
\end{equation}
Thus, $lrd(n,\alpha_{n})/|I_{\alpha_n}^{n}|+|I_{\alpha_n+1}^{n}|$
converges to a value $d_P$,
when $n$ tends to infinity.
Let us prove that $\psi$ is differentiable at $P$ and that
$d\psi(P)=2d_P$. Let $L=[x,y]$ be any   interval
such that the point $P \in L$.
Take the largest $m \ge 1$ such that
there is a grid interval $I_\gamma^m$
with the property that
$L \subset   \cup_{j=-1,0,1} I_{\gamma+j}^m$.
By the  bounded geometry property of a grid, there is $C \ge 1$,
not depending upon $P$, $L$ and $I_\gamma^m$,
such that
\begin{equation}
\label{BAC69}
C^{-1} \le
\frac{|I_\gamma^m|}{|L|}
\le C \ .
\end{equation}
For every $n \ge m$, take the smallest sequence of adjacent
grid intervals
$I_{\beta_n}^n, \ldots, I_{\beta_n+i_n}^n$,
at level $n$,
such that $L \subset \cup_{i=0}^{i_n} I_{\beta_n+i}^n \subset
\cup_{j=-1,0,1} I_{\gamma+j}^m$. Hence, by definition of the
logarithmic ratio distortion, we get
$$
\psi(x)=\lim_{n \to \infty} ldh(I_{\beta_n}^n)
$$
and
$$
\psi(y)=\lim_{n \to \infty} ldh(I_{\beta_n+i_n}^n) \ .
$$
Therefore,
\begin{eqnarray}
\label{BAC9}
\frac{\psi(y)-\psi(x)}{y-x} & = &
\lim_{n \to \infty} \frac{ldh(I_{\beta_n+i_n}^n)-ldh(I_{\beta_n}^n)}{y-x}
\nonumber\\
& = &
\lim_{n \to \infty} \frac{\sum_{i=0}^{i_n-1}lrd(I_{\beta_n+i}^n)}{y-x} \ .
\end{eqnarray}
By inequalities (\ref{BAC7}) and (\ref{BAC69}), for every
$I_{\beta_n+i}^n \subset I_{\gamma+j(i)}^m$ ,we get
\begin{eqnarray}
\label{BAC11}
lrd(n,\beta_{n}+i)   &  \in   &
\left(|I_{\beta_n+i}^{n}|+|I_{\beta_n+i+1}^{n}|\right)
\left(
\frac{lrd(m,\gamma+j(i))}{|I_{\gamma+j(i)}^{m}|+|I_{\gamma+j(i)+1}^{m}|}
\pm \cO(|I_{\gamma+j(i)}^m|^{r}) \right) \nonumber \\
& \subset & \left(|I_{\beta_n+i}^{n}|+|I_{\beta_n+i+1}^{n}|\right)
\left( d_P  \pm \cO(|L|^{r}) \right) \ .
\end{eqnarray}
Putting together (\ref{BAC9}) and (\ref{BAC11}),
we obtain that
\begin{eqnarray}
\label{BA33}
\frac{\psi(y)-\psi(x)}{y-x} & \in &
\lim_{n \to \infty} ( d_P  \pm \cO(|L|^{r})
\frac{\sum_{i=0}^{i_n-1}
|I_{\beta_n+i}^{n}|+|I_{\beta_n+i+1}^{n}|}{y-x} \nonumber \\
& \subset &
\lim_{n \to \infty} ( d_P  \pm \cO(|L|^{r})
\frac{|I_{\beta_n}^{n}|+|I_{\beta_n+i_n}^{n}|+
2\sum_{i=1}^{i_n-1}|I_{\beta_n+i}^{n}|}{y-x}
\nonumber \\
& \subset &
2  d_P  \pm \cO(|L|^{r}) \ .
\end{eqnarray}
Therefore, for every $P \in I$, the homeomorphism $\psi$ is
differentiable at $P$ and $d\psi(P)=2d_P$. Let us check that $d\psi$
is $r$-H\"older continuous. For every $P,P' \in I$, let $L$ be the
closed interval $[P,P']$. Using (\ref{BA33}), we obtain that
\begin{eqnarray*}
\label{BAC73}
d\psi(P')-d\psi(P) &\in &
\frac{\psi(P')-\psi(P)}{P'-P} - \frac{\psi(P')-\psi(P)}{P'-P}
\pm
\cO(|L|^{r}) \\
&\subset & \pm \cO(|L|^{r}) \ ,
\end{eqnarray*}
and so $d\psi$ is $r$-H\"older continuous.

\bigskip
\noindent
{\it Case $r=1$}: By the  above argument,  $h$ is $C^{2+s}$ for every
$0<s<1$ and so, in particular, $h$ is $C^{1+Lipschitz}$. Thus, by
Lemma \ref{r1+}, for every $n \ge 1$ and
every $1 \le \beta \le \Omega(n)-1$ we get that
$$
|lrd(n,\beta)| \le \cO(|I_{\beta}^n|) \ ,
$$
which implies that  inequality  (\ref{batatas53}) is also satisfied   for
$r=1$.
Now, a similar argument to the one above gives that
$d\psi$ is Lipschitz.
\qed

\subsection{Proof of Theorem \ref{th13131313}}
\label{pppttt}

In this section, we prove Theorem \ref{th13131313}.

\bigskip

\noindent {\bf Proof of Theorem \ref{th13131313}:} The equivalences
presented for quasisymmetric homeomorphisms follow from  Lemma
\ref{lll0} with respect to ratio distortion and from Lemma
\ref{lll1}  with respect to cross ratio distortion, noting that the
ratios $r(n,\beta)$ and the cross ratios $cr(n,\beta)$ are uniformly
bounded by the
bounded geometry property of the grid.
The equivalences presented for uniformly asymptotically affine
(uaa) homeomorphisms follow from  Lemma  \ref{gfsdgfgd33} with
respect to ratio distortion and from Lemma  \ref{lll111}  with
respect to cross ratio distortion. The equivalences presented for
$C^{1+\alpha}$, $C^{1+Lipschitz}$ and affine diffeomorphisms
follow from Lemma \ref{r1+} with respect to ratio distortion  and
from Lemma \ref{lll11179}
with respect to cross ratio distortion.
% The equivalences presented for $C^{1+Zigmund}$ and $C^{1+zigmund}$
%diffeomorphisms    follow from Lemma \ref{lll111Z}.
The equivalences presented for  $C^{2+\alpha}$ and $C^{2+Lipschitz}$
diffeomorphisms    follow from Lemma \ref{lll111r2}. \qed

\par\bigskip
\par\bigskip
\noindent {\Large\bf Acknowledgements} We would like to thank
David Rand, Nils Tongring and Fl\'avio Ferreira all the useful discussions.
A.
Pinto would like to thank  CUNY, IHES, IMPA, University of Warwick
and SUNY for their hospitality, and
   Calouste Gulbenkian Foundation, PRODYN-ESF,   FCT of
  MCT, and CMUP for their
financial support.

\end{document}

%% file: imsmark.tex
\def\IMSmarkvadjust{0 pt}
\def\IMSmarkhadjust{0 pt}
\def\IMSmarkhpadding{0 pt}
\def\IMSpubltext{Published in modified form:}
\def\SBIMSMark#1#2#3{
 \font\SBF=cmss10 at 10 true pt
 \font\SBI=cmssi10 at 10 true pt
 \setbox0=\hbox{\SBF \hbox to \IMSmarkhpadding{\relax}
                Stony Brook IMS Preprint \##1}
 \setbox2=\hbox to \wd0{\hfil \SBI #2}
 \setbox4=\hbox to \wd0{\hfil \SBI #3}
 \setbox6=\hbox to \wd0{\hss
             \vbox{\hsize=\wd0 \parskip=0pt \baselineskip=10 true pt
                   \copy0 \break%
                   \copy2 \break% 
                   \copy4 \break}}
 \dimen0=\ht6   \advance\dimen0 by \vsize \advance\dimen0 by 8 true pt
                \advance\dimen0 by -\pagetotal
	        \advance\dimen0 by \IMSmarkvadjust
 \dimen2=\hsize \advance\dimen2 by .25 true in
	        \advance\dimen2 by \IMSmarkhadjust

%
%   Check for publication info
%
%  \newread\jref
  \openin2=publishd.tex
  \ifeof2\setbox0=\hbox to 0pt{}
  \else 
     \setbox0=\hbox to 3.1 true in{
                \vbox to \ht6{\hsize=3 true in \parskip=0pt  \noindent  
                {\SBI \IMSpubltext}\hfil\break
                \input publishd.tex 
                \vfill}}
  \fi
  \closein2
  \ht0=0pt \dp0=0pt
 \ht6=0pt \dp6=0pt
 \setbox8=\vbox to \dimen0{\vfill \hbox to \dimen2{\copy0 \hss \copy6}}
 \ht8=0pt \dp8=0pt \wd8=0pt
 \copy8
 \message{*** Stony Brook IMS Preprint #1, #2. #3 ***}
}